\documentclass[preprint,12pt]{elsarticle}

\usepackage{amsthm}
\usepackage{amssymb,amsfonts}
\usepackage[tbtags]{amsmath} %tbtags：用于将split环境下的公式编号由居中放到最后一行公式后，fleqn：用于将gather环境里的公式由居中对齐改为左对齐
\usepackage{booktabs}
\allowdisplaybreaks  %使用amsmath后，公式环境不能自动分页，用此命令即可

\usepackage{caption}
\newtheorem{thm}{Theorem}
\newtheorem{lem}{Lemma}
\newtheorem{pro}{Proposition}

\newdefinition{rmk}{Remark}
\newproof{pf}{Proof}
\newproof{pot}{Proof of Theorem \ref{thm2}}

\begin{document}

\begin{frontmatter}

\title{Heavy-traffic Asymptotics of Priority Polling System with Threshold Service Policy}

\author{Liu Zaiming}
\ead{math\_lzm@csu.edu.cn}
\author{Chu Yuqing}
\ead{chuyuqing@csu.edu.cn}
\author{Wu Jinbiao}
\ead{Corresponding author: wujinbiao@csu.edu.cn}

\address{Department of Mathematics and Statistics, Central
South University, Changsha, Hunan 410083, PR China}

\begin{abstract}
In this paper, by the singular-perturbation technique, we
investigate the heavy-traffic behavior of a priority polling system
consisting of three $M/M/1$ queues with threshold policy. It turns
out that the scaled queue-length of the critically loaded queue is
exponentially distributed, independent of that of the stable queues.
In addition, the queue lengths of stable queues possess the same
distributions as a priority polling system with $N$-policy vacation.
Based on this fact, we provide the exact tail asymptotics of the
vacation polling system to approximate the tail distribution of the
queue lengths of the stable queues, which shows that it has the same
prefactors and decay rates as the classical $M/M/1$ preemptive
priority queues. Finally, a stochastic simulation is taken to test
the results aforementioned.
\end{abstract}


\end{frontmatter}

\section{Introduction}\label{sec1}
The
study of the two-queue priority polling system is motivated
by its wide applications in computer and communication systems, such as ATM (Asynchronous Transfer Mode) switch systems and
network standards like DQDB(Distributed Queue Dual Bus). ATM
involves two different types of traffic: real time traffic(voice,
video) and non-real time traffic (data), which also
need different types of QoS (Quality of Service) standard. By setting the threshold
parameter, a higher priority is offered
to real time traffic to shorten its delay and the delay of non-real
time traffic is kept in a valid regime, which turns out to be a flexible way to control the operation of the whole system.

Lee
and Sengupta first investigated the threshold-based priority systems in \cite{Lee1993}. Later, a special case of two-queue $M/M/1$ polling
system with threshold policy was studied by Boxma, Koole and Mitrani in
\cite{Boxma1995,Boxma1995a}. The model was further
extended with switch-over times by Deng et al.
 in \cite{Deng2001,Deng2001a} and with one more server by Feng in
\cite{Feng2001}.

In \cite{Liu2014}, we concerned with a
three-queue model under threshold policy. The motivation stems from
\cite{Landry1993}, in which Landry and Stavrakakis proposed a third
type of traffic so called control traffic with Head-of-Line (HoL) in the integrated ATM environment, which involves
critical network control and reservation information. In this paper, we focus on the heavy-traffic limits when there is a single critically loaded queue.

Using the singular-perturbation technique, we derive the lowest-order asymptotic of the joint queue-length distribution in terms of a small positive parameter measuring the closeness of the system to instability.
The singular-perturbation technique was first applied to investigate the heavy-traffic behavior of interacting queues in \cite{Morrison2010}. Later, Boon and Winands \cite{Boon2013} used this technique to a model with $k$-limited policies and presented the heavy-traffic behavior.
It is noted that the singular-perturbation technique can be easily extended to a multi-queue system since it only needs the balance equations.

With the singular-perturbation technique, we conclude that the queue lengths in the stable queues have the same joint distribution as Model {\uppercase\expandafter{\romannumeral2}}, a preemptive priority polling system with $N$-policy vacation. In general, no closed-form expressions for the steady-state probabilities in Model {\uppercase\expandafter{\romannumeral2}} can be obtained. Using the Kernel method, which is reported detailedly in \cite{Li2009,Li2011}, we present the exact tail asymptotics of queue lengths in Model {\uppercase\expandafter{\romannumeral2}}, which can further approximate the tail asymptotics of the stable queues.

The remainder of this paper is organized as follows.
In Section \ref{sec2}, the model and some notations are introduced.
In Section \ref{sec3}, the singular-perturbation technique is applied to derive the heavy-traffic limits and the detailed derivation is carried out in Section \ref{sec4}.
In Section \ref{sec5}, we provide the exact tail asymptotics of queue lengths in Model {\uppercase\expandafter{\romannumeral2}} to approximate the tail asymptotics of the stable queues.
In Section \ref{sec6}, a simulation is undertaken to evaluate the heavy-traffic asymptotics.
We finally conclude the whole procedure and propose some topics for further research in Section \ref{sec7}.

 \section{Model Description}\label{sec2}
 We consider a polling model with single server consisting of
 three queues $Q_{1}$, $Q_{2}$, $Q_{3}$. We refer to the customers queueing in $Q_{i}$ as the type $i$ customers, $i=1,2,3$. The buffer capacity of each queue is infinite.
 Customers arrive at $Q_{i}$ independently according
 to a Poisson process with rate $\lambda_{i}$. For type $i$ customers, the service times
 are mutually independent and all follow an exponential distribution with rate $\mu_{i}$.
 $Q_{1}$ has
 the HoL priority and $Q_{2}$ has a higher priority over $Q_{3}$. In each queue customers are served according to FCFS discipline. We assume that
 the arrival processes and the service processes are
independent.
 The service discipline is described as follows.
 \begin{enumerate}[1.]
    \item $Q_{1}$ is served exhaustively, which means that the server serves the customers in $Q_{1}$
    until it is empty and then switches to $Q_{2}$;
    \item When the server is serving a customer in $Q_{2}$, if a type 1 customer arrives,
    then the server switches to $Q_{1}$ immediately, otherwise,
    it continues serving the customers in $Q_{2}$ until $Q_{2}$ becomes empty and then switches to $Q_{3}$;
    \item When the server is serving a customer in $Q_{3}$,
    if a type 1 customer arrives, then the server switches to $Q_{1}$
    immediately,
    if the size of $Q_{2}$ reaches a given threshold $N$ and $Q_1$ is empty,
    then the server switches to $Q_{2}$ immediately, otherwise,
    it continues serving the customers in $Q_{3}$ until $Q_{3}$ becomes empty and
    then switches to $Q_{2}$.
\end{enumerate}
It is assumed that all the switches are instantaneous.
 In addition, the switches caused by the threshold push the customer
 undergoing service to the head of the queue and the service of the interrupted customer
resumes from the beginning.

The traffic load of $Q_{i}$ is denoted
by $\rho_i = \lambda_{i}/\mu_{i}$, $i=1,2,3$.
We assume the ergodicity condition of the system $\rho=\rho_1+\rho_2+\rho_3<1$ is satisfied.

Let $X_i(t)$ be the number of customers in $Q_i$ at time $t$,
and $S(t)$ be the position of the server at time $t$ with $S(t)\in \{1,2,3\}$.
The associated stochastic process $\{Y(t),t\geq 0\}=\left\{\big(X_1(t),X_2(t),X_3(t),S(t)\big),t\geq 0\right\}$ is an aperiodic and irreducible four-dimensional Markov process.
Let $X_i$ ($i=1,2,3$) be the steady-state queue length of $Q_i$ and $S$ be the steady-state position of the server.
Define the stationary probabilities:
\begin{equation*}
% \nonumber to remove numbering (before each equation)
p_{s}(x_1,x_2,x_3) = \lim_{t\rightarrow \infty}Pr\{Y(t)=(x_1,x_2,x_3,s)\}, \ \ \  s=1,2,3.
\end{equation*}

We study the heavy-traffic limits of the joint queue-length distribution by increasing the arrival rate $\lambda_3$ so as to  $\rho\rightarrow 1^{-}$, while keeping $\lambda_1\neq 0$, $\lambda_2\neq 0$ and $\mu_1, \mu_2, \mu_3$ fixed. When $\rho\rightarrow 1^{-}$, $Q_3$ becomes critically loaded, whereas $Q_1$ and $Q_2$ remain stable since $Q_1$ and $Q_2$ have higher priorities over $Q_3$.

The single-perturbation technique is implemented here.
We first apply a perturbation to $\lambda_3$ in the balance equations, in which case $Q_3$ is close to becoming critically loaded.
Then we solve the lowest order terms in the balance equations to obtain the queue-length distributions of the stable queues $Q_1$ and $Q_2$. At last we solve the first-order and second-order terms to get a differential equation and compute the scaled number of customers in $Q_3$.

Applying the Markov property, we obtain the following balance equations when $x_3\geq 2$:
\begin{gather}
    \begin{split}\label{eq:21}
      &(\lambda_1+\lambda_2 +\lambda_3+\mu_1)p_1(x_1,x_2,x_3) \\
                           &=\lambda_1 p_1(x_1-1,x_2,x_3)\delta(x_1\geq 2)+\lambda_2 p_1(x_1,x_2-1,x_3)\delta(x_2\geq 1)\\
                            &\ \ \ +\lambda_3 p_1(x_1,x_2,x_3-1)+\lambda_1 p_2(0,x_2,x_3)\delta(x_1=1,x_2\geq 1)\\
                            &\ \ \ +\mu_1 p_1(x_1+1,x_2,x_3)+\lambda_1 p_3(0,x_2,x_3)\delta(x_1=1,x_2<N), x_1\geq 1, x_2\geq 0,
    \end{split}\\
    \begin{split}\label{eq:22}
      &(\lambda_1+\lambda_2 +\lambda_3+\mu_2)p_2(0,x_2,x_3) \\
                           &=\lambda_2 p_2(0,x_2-1,x_3)\delta(x_2\geq 2)+\lambda_3 p_2(0,x_2,x_3-1)+\mu_1 p_1(1,x_2,x_3)\\
                           &\ \ \ +\lambda_2 p_3(0,N-1,x_3)\delta(x_2=N)+\mu_2 p_2(0,x_2+1,x_3), \ \ \ x_2\geq 1,
    \end{split}\\
    \begin{split}\label{eq:23}
      &(\lambda_1+\lambda_2+\lambda_3+\mu_3)p_3(0,x_2,x_3) \\
                           & =\lambda_2 p_3(0,x_2-1,x_3)\delta(x_2\geq 1)+\lambda_3 p_3(0,x_2,x_3-1)+\mu_3 p_3(0,x_2,x_3+1)\\
                           &\ \ \ +\left[\mu_1 p_1(1,0,x_3)+\mu_2 p_2(0,1,x_3)\right]\delta(x_2=0), \ \ \ 0\leq x_2 \leq N-1,
    \end{split}
\end{gather}
where $\delta(\cdot)$ is Kronecker function.

In the above equations, we have omitted the parts for $x_3=0$ and $x_3=1$ which do not play a role after the perturbation since $X_3$ tends to infinity as $Q_3$ becomes critically loaded and the probability of $Q_3$ being empty or 1 goes to zero.

Throughout the paper, we adopt the standard notations: a function $F(x)$ is $o(x)$ if $F(x)/x\rightarrow 0$ as $x\rightarrow 0$; a function $F(x)$ is $\mathcal{O}(x)$ if there exists a $c\geq0 $ such that $F(x)/x\rightarrow c$ as $x\rightarrow 0$ while $\mathcal{O}(1)$ is a constant time complexity; functions $f(n)$ and $g(n)$ of nonnegative integers $n$, $f(n)\sim g(n)$ means $\lim_{n\rightarrow \infty}\frac{f(n)}{g(n)}=1$.

\section{Perturbation}\label{sec3}
From the stability condition the system becomes unstable as $\rho_3\rightarrow 1-\rho_1-\rho_2$, i.e. $\lambda_3\rightarrow \mu_3(1-\rho_1-\rho_2)$. Therefore it is assumed that
\begin{equation}\label{eq:31}
\lambda_3=\mu_3(1-\rho_1-\rho_2)-\varepsilon\omega,\ \ \ \omega>0, 0<\varepsilon \ll 1.
\end{equation}

Let $\zeta=\varepsilon x_3$, and
\begin{equation}\label{eq:32}
  p_{s}(x_1,x_2,x_3)=p_{s}(x_1,x_2,\zeta/\varepsilon)=\varepsilon \phi_{s,(x_1,x_2)}(\zeta,\varepsilon),\ \ \ 0<\zeta=\mathcal{O}(1),s=1,2,3.
\end{equation}

Taking \eqref{eq:31} and \eqref{eq:32} into the balance equations \eqref{eq:21}-\eqref{eq:23} and then taking the Taylor expansion, we obtain
\begin{gather}
    \begin{split}\label{eq:36}
      &(\lambda_1+\lambda_2 +\mu_1)\phi_{1,(x_1,x_2)}(\zeta,\varepsilon) \\
                           &=\lambda_1 \phi_{1,(x_1-1,x_2)}(\zeta,\varepsilon)\delta(x_1\geq 2)+\lambda_2 \phi_{1,(x_1,x_2-1)}(\zeta,\varepsilon)\delta(x_2\geq 1)\\
                            &\ \ \ -(\mu_3(1-\rho_1-\rho_2)-\varepsilon\omega) \left(\varepsilon\frac{\partial\phi_{1,(x_1,x_2)}(\zeta,\varepsilon)}{\partial\zeta}-\frac{\varepsilon^2}{2}\frac{\partial^2\phi_{1,(x_1,x_2)}(\zeta,\varepsilon)}{\partial\zeta^2}\right)\\
                            &\ \ \ +\lambda_1 \phi_{2,(0,x_2)}(\zeta,\varepsilon)\delta(x_1=1,x_2\geq 1)+\mu_1 \phi_{1,(x_1+1,x_2)}(\zeta,\varepsilon)\\
                            &\ \ \ +\lambda_1 \phi_{3,(0,x_2)}(\zeta,\varepsilon)\delta(x_1=1,x_2<N)+o(\varepsilon^2),\ \ \  x_1\geq 1, x_2\geq 0,
    \end{split}\\
    \begin{split}\label{eq:37}
      &(\lambda_1+\lambda_2 +\mu_2)\phi_{2,(0,x_2)}(\zeta,\varepsilon) \\
                           & =(\mu_3(1-\rho_1-\rho_2)-\varepsilon\omega) \left(-\varepsilon\frac{\partial\phi_{2,(0,x_2)}(\zeta,\varepsilon)}{\partial\zeta}+\frac{\varepsilon^2}{2}\frac{\partial^2\phi_{2,(0,x_2)}(\zeta,\varepsilon)}{\partial\zeta^2}\right)\\
                           &\ \ \ +\lambda_2 \phi_{2,(0,x_2-1)}(\zeta,\varepsilon)\delta(x_2\geq 2)+\mu_1 \phi_{1,(1,x_2)}(\zeta,\varepsilon)\\
                           &\ \ \ +\lambda_2 \phi_{3,(0,N-1)}(\zeta,\varepsilon)\delta(x_2=N)+\mu_2 \phi_{2,(0,x_2+1)}(\zeta,\varepsilon)+o(\varepsilon^2), \ \ x_2\geq 1,
    \end{split}\\
    \begin{split}\label{eq:38}
      &(\lambda_1+\lambda_2)\phi_{3,(0,x_2)}(\zeta,\varepsilon) \\
                           & =(\mu_3(1-\rho_1-\rho_2)-\varepsilon\omega) \left(-\varepsilon\frac{\partial\phi_{3,(0,x_2)}(\zeta,\varepsilon)}{\partial\zeta}+\frac{\varepsilon^2}{2}\frac{\partial^2\phi_{3,(0,x_2)}(\zeta,\varepsilon)}{\partial\zeta^2}\right)\\
                           &\ \ \ +\mu_3\left(\varepsilon\frac{\partial\phi_{3,(0,x_2)}(\zeta,\varepsilon)}{\partial\zeta}+\frac{\varepsilon^2}{2}\frac{\partial^2\phi_{3,(0,x_2)}(\zeta,\varepsilon)}{\partial\zeta^2}\right)\\
                           &\ \ \ +\left[\mu_1 \phi_{1,(1,0)}(\zeta,\varepsilon)+\mu_2 \phi_{2,(0,1)}(\zeta,\varepsilon)\right]\delta(x_2=0)\\
                           &\ \ \ +\lambda_2 \phi_{3,(0,x_2-1)}(\zeta,\varepsilon)\delta(x_2\geq 1)+o(\varepsilon^2), \ \ \ 0\leq x_2 \leq N-1.
    \end{split}
\end{gather}

It is noted that $\lambda_3$ only plays a role in equations for $\mathcal{O}(\varepsilon)$ terms and higher. Throughout the paper, we do Taylor expansions of $\phi_{s,(x_1,x_2)}(\zeta,\varepsilon)$ ($s=1,2,3$) in powers of $\varepsilon$ as follows
\begin{equation}\label{eq:380}
 \phi_{s,(x_1,x_2)}(\zeta,\varepsilon)=\phi_{s,(x_1,x_2)}^{(0)}(\zeta)+\varepsilon\phi_{s,(x_1,x_2)}^{(1)}(\zeta)+o(\varepsilon^2), \ \ \ s=1,2,3.
\end{equation}
In the next section the lowest order terms of the resulting equations after Taylor expansions are equated to find expressions for $\phi_{s,(x_1,x_2)}^{(0)}(\zeta)$ ($s=1,2,3$), subsequently the first-order and second-order terms are equated to find the scaled queue-length distribution of $Q_3$.

For convenience, we introduce the corresponding probability generating functions(PGFs):
\begin{align*}
 &Q_{1}^{(j)}(x,y,\zeta)=\sum_{x_1=1}^{\infty}\sum_{x_2=0}^{\infty}\phi_{1,(x_1,x_2)}^{(j)}(\zeta)x^{x_1-1}y^{x_2},\ \ \ j=0,1,\\ &Q_{2}^{(j)}(y,\zeta)=\sum_{x_2=1}^{\infty}\phi_{2,(0,x_2)}^{(j)}(\zeta)y^{x_2-1}, \ \ \ j=0,1,\\
 &Q_{3}^{(j)}(y,\zeta)=\sum_{x_2=0}^{N-1}\phi_{3,(0,x_2)}^{(j)}(\zeta)y^{x_2},\ \ \  j=0,1,\\
&Q_{1}(x,y,\zeta,\varepsilon)=\sum_{x_1=1}^{\infty}\sum_{x_2=0}^{\infty}\phi_{1,(x_1,x_2)}(\zeta,\varepsilon)x^{x_1-1}y^{x_2},\\ &Q_{2}(y,\zeta,\varepsilon)=\sum_{x_2=1}^{\infty}\phi_{2,(0,x_2)}(\zeta,\varepsilon)y^{x_2-1},\\
 &Q_{3}(y,\zeta,\varepsilon)=\sum_{x_2=0}^{N-1}\phi_{3,(0,x_2)}(\zeta,\varepsilon)y^{x_2}.
\end{align*}

\section{Model analysis}\label{sec4}
\subsection{Equating the lowest-order terms}\label{subsec4.1}
Equating the lowest-order terms of the resulting equation after the Taylor expansions of \eqref{eq:36}-\eqref{eq:38}, we obtain
\begin{gather}
    \begin{split}\label{eq:39}
      &(\lambda_1+\lambda_2 +\mu_1)\phi_{1,(x_1,x_2)}^{(0)}(\zeta) \\
                           &=\lambda_1 \phi_{1,(x_1-1,x_2)}^{(0)}(\zeta)\delta(x_1\geq 2)+\lambda_2 \phi_{1,(x_1,x_2-1)}^{(0)}(\zeta)\delta(x_2\geq 1)\\
                            &\ \ \ +\lambda_1 \phi_{2,(0,x_2)}^{(0)}(\zeta)\delta(x_1=1,x_2\geq 1)+\mu_1 \phi_{1,(x_1+1,x_2)}^{(0)}(\zeta)\\
                            &\ \ \ +\lambda_1 \phi_{3,(0,x_2)}^{(0)}(\zeta)\delta(x_1=1,x_2<N),\ \ \  x_1\geq 1, x_2\geq 0,
    \end{split}\\
    \begin{split}\label{eq:310}
      &(\lambda_1+\lambda_2 +\mu_2)\phi_{2,(0,x_2)}^{(0)}(\zeta) \\
                           &=\lambda_2 \phi_{2,(0,x_2-1)}^{(0)}(\zeta)\delta(x_2\geq 2)+\mu_1 \phi_{1,(1,x_2)}^{(0)}(\zeta)\\
                           &\ \ \ +\lambda_2 \phi_{3,(0,N-1)}^{(0)}(\zeta)\delta(x_2=N)+\mu_2 \phi_{2,(0,x_2+1)}^{(0)}(\zeta),\ \ \  x_2\geq 1,
    \end{split}\\
    \begin{split}\label{eq:311}
      (\lambda_1+\lambda_2)\phi_{3,(0,x_2)}^{(0)}(\zeta)=\lambda_2 \phi_{3,(0,x_2-1)}^{(0)}(\zeta), \ \ \ 1\leq x_2 \leq N-1,
    \end{split}\\
\begin{split}\label{eq:312}
      (\lambda_1+\lambda_2)\phi_{3,(0,0)}^{(0)}(\zeta)=\mu_1 \phi_{1,(1,0)}^{(0)}(\zeta)+\mu_2 \phi_{2,(0,1)}^{(0)}(\zeta).
    \end{split}
\end{gather}

We introduce $P_0(\zeta)$ and $\pi_{s,(x_1,x_2)}^{(0)}$ such that
\begin{gather*}
\phi_{s,(x_1,x_2)}^{(0)}(\zeta)=\pi_{s,(x_1,x_2)}^{(0)}P_0(\zeta),\qquad s=1,2,3, \label{eq:313}\\ \sum_{x_1=1}^{\infty}\sum_{x_2=0}^{\infty}\pi_{1,(x_1,x_2)}^{(0)}+\sum_{x_2=1}^{\infty}\pi_{2,(0,x_2)}^{(0)}+\sum_{x_2=0}^{N-1}\pi_{3,(0,x_2)}^{(0)}=1.\label{eq:314}
\end{gather*}
Define
\begin{align*}
 &L_{1}^{(0)}(x,y)=\sum_{x_1=1}^{\infty}\sum_{x_2=0}^{\infty}\pi_{1,(x_1,x_2)}^{(0)}x^{x_1-1}y^{x_2},
 &&L_{2}^{(0)}(y)=\sum_{x_2=1}^{\infty}\pi_{2,(0,x_2)}^{(0)}y^{x_2-1},\\
 &L_{3}^{(0)}(y)=\sum_{x_2=0}^{N-1}\pi_{3,(0,x_2)}^{(0)}y^{x_2}. & &
\end{align*}
Then it is clear that
\begin{align}
&Q_{1}^{(0)}(x,y,\zeta)=L_1^{(0)}(x,y)P_0(\zeta),\label{eq:3181}\\
&Q_{2}^{(0)}(y,\zeta)=L_2^{(0)}(y)P_0(\zeta),\label{eq:3191}\\
&Q_{3}^{(0)}(y,\zeta)=L_{3}^{(0)}(y)P_0(\zeta).\label{eq:3201}
\end{align}
From \eqref{eq:311}, we get
\begin{equation}\label{eq:a2}
  L_{3}^{(0)}(y)=\sum_{x_2=0}^{N-1}(r_2y)^{x_2}\pi_{3,(0,0)}^{(0)}=H(y)\pi_{3,(0,0)}^{(0)}=\beta(y)L_3^{(0)}(1),
\end{equation}
where $r_2=\frac{\lambda_2}{\lambda_1+\lambda_2}$, $H(y)=\frac{(r_2y)^N-1}{r_2y-1}$ and $\beta(y)=\frac{H(y)}{H(1)}$.

Using the PGFs to rewrite the balance equations \eqref{eq:39} and \eqref{eq:310} leads to
\begin{align}
      &xK(x,y)L_1^{(0)}(x,y) =\lambda_1 x[yL_2^{(0)}(y)+L_3^{(0)}(y)]-\mu_1 L_1^{(0)}(0,y),\label{eq:a5}\\
     &ya(y)L_2^{(0)}(y)=\mu_1 L_1^{(0)}(0,y)-[\lambda_1+\lambda_2(1-y)]L_3^{(0)}(y),\label{eq:a6}
\end{align}
where
\begin{align*}
&K(x,y)=\lambda_1(1-x)+\lambda_2(1-y)+\mu_1\left(1-\frac{1}{x}\right),\\
 &a(y)=\lambda_1+\lambda_2(1-y)+\mu_2\left(1-\frac{1}{y}\right).
\end{align*}

 Clearly, for every $|y|\leq 1$, the kernel $xK(x,y)$ has a unique zero: $x=\alpha(y)$. Applying the Kernel method to \eqref{eq:a5} and \eqref{eq:a6}, it is easy to get
 \begin{align}
&L_1^{(0)}(x,y)=\frac{\lambda_1\mu_2[x-\alpha(y)](y-1)}{xK(x,y)[ya(y)-\lambda_1y\alpha(y)]}L_3^{(0)}(y),\label{eq:318}\\
&L_2^{(0)}(y)=\frac{\lambda_1(\alpha(y)-1)+\lambda_2(y-1)}{ya(y)-\lambda_1y\alpha(y)}L_3^{(0)}(y),\label{eq:319}
\end{align}
Letting $y\rightarrow 1$ and then letting $x\rightarrow 1$ in \eqref{eq:318} and \eqref{eq:319}, with L'H\^{o}pital's rule, we obtain $L_1^{(0)}(1,1)=\frac{\rho_1}{1-\rho_1-\rho_2}L_3^{(0)}(1)$ and $L_2^{(0)}(1)=\frac{\rho_2}{1-\rho_1-\rho_2}L_3^{(0)}(1)$. By the normalizing condition, it is easy to get $L_3^{(0)}(1)=1-\rho_1-\rho_2$. Therefore, we have $L_1^{(0)}(1,1)=\rho_1$ and $L_2^{(0)}(1)=\rho_2$. Moreover,
\begin{equation}\label{eq:320}
L_{3}^{(0)}(y)=\beta(y)(1-\rho_1-\rho_2).
\end{equation}

It is not hard to see that equations \eqref{eq:318}-\eqref{eq:320} actually state an $M/M/1$ preemptive priority polling system with $N$-policy vacation, denoted as Model \uppercase\expandafter{\romannumeral2} for short, described as follows:

 There are two classes of customers in the system, the high- and low-priority customers, arriving independently according to two Poisson processes with rates $\lambda_1$ and $\lambda_2$, respectively.
 Each class of customer is served according to the FCFS discipline.
 The server takes a vacation once the system empties and goes back to work once the size of the low-priority customers reaches $N$ or there is a high-priority customer's arrival. The high-priority customers have preemptive priorities over the low-priority customers just like in the classical two-queue preemptive priority queueing system.
 Both classes of customers require an exponential amount of service times and are served with service rates $\mu_1$ and $\mu_2$, respectively.
 All service times are independent and also independent of the arrival processes.

We determine the unkown expression of $P_0(\zeta)$ in the rest of this section.

\subsection{Equating the first-order terms}\label{subsec4.2}
In this subsection, by equating the first-order terms of the resulting equations after the Taylor expansion of the perturbed balance equations \eqref{eq:36}-\eqref{eq:38}, we present an equation in Proposition \ref{pro1}.

\begin{pro}\label{pro1}
\begin{multline*}
(1-\rho_1-\rho_2)\left[Q_{1}^{(1)}(1,1,\zeta)+Q_{2}^{(1)}(1,\zeta)+Q_{3}^{(1)}(1,\zeta)\right]-Q_{3}^{(1)}(1,\zeta)\\
=-\left[\frac{\mu_3}{\mu_1}\rho_1+\frac{\mu_3}{\mu_2}\rho_2\right]P_0^{'}(\zeta).
\end{multline*}
\end{pro}

\begin{pf}
Taking the PGF of the first-order terms of the resulting equations after the Taylor expansion of \eqref{eq:36}-\eqref{eq:38}, we have
\begin{gather}
    \begin{split}\label{eq:3f5}
      xK(x,y)Q_{1}^{(1)}(x,y,\zeta)=&\lambda_1xy Q_{2}^{(1)}(y,\zeta)-\mu_1Q_{1}^{(1)}(0,y,\zeta)+\lambda_1xQ_{3}^{(1)}(y,\zeta) \\
                           &-\mu_3x(1-\rho_1-\rho_2)L_{1}^{(0)}(x,y)P_0^{'}(\zeta),
    \end{split}\\
    \begin{split}\label{eq:3f6}
       ya(y)Q_{2}^{(1)}(y,\zeta)=&\lambda_2y^N\phi_{3,(0,x_2-1)}^{(1)}(\zeta)-\mu_1Q_{1}^{(1)}(0,0,\zeta)-\mu_2Q_{2}^{(1)}(0,\zeta)\\
                            &+\mu_1Q_{1}^{(1)}(0,y,\zeta)-\mu_3y(1-\rho_1-\rho_2)L_{2}^{(0)}(y)P_0^{'}(\zeta),
    \end{split}\\
    \begin{split}\label{eq:3f7}
      [\lambda_1+\lambda_2(1-y)]Q_{3}^{(1)}(y,\zeta)=&-\lambda_2y^N\phi_{3,(0,x_2-1)}^{(1)}(\zeta)+\mu_1Q_{1}^{(1)}(0,0,\zeta)\\
      &+\mu_2Q_{2}^{(1)}(0,\zeta)+\mu_3(\rho_1+\rho_2)L_{3}^{(0)}(y)P_0^{'}(\zeta).
    \end{split}
\end{gather}
Applying the Kernel method to \eqref{eq:3f5}-\eqref{eq:3f7}, after some elementary calculations, we get
\begin{gather}
\begin{split}
Q_{1}^{(1)}(x,y,\zeta)=&\frac{\lambda_1[x-\alpha(y)]}{xK(x,y)}[yQ_{2}^{(1)}(y,\zeta)+Q_{3}^{(1)}(y,\zeta)]\\
&-\frac{\mu_3(1-\rho_1-\rho_2)}{xK(x,y)}[xL_{1}^{(0)}(x,y)-\alpha(y)L_{1}^{(0)}(\alpha(y),y)]P_0^{'}(\zeta).
\end{split}\label{eq:3f9}\\
\begin{split}\label{eq:3f11}
 &[ya(y)-\lambda_1\alpha(y)y]Q_{2}^{(1)}(y,\zeta)+[\lambda_1(1-\alpha(y))+\lambda_2(1-y)]Q_{3}^{(1)}(y,\zeta)\\
    &=\mu_3\left\{(\rho_1+\rho_2)L_{3}^{(0)}(y)-(1-\rho_1-\rho_2)\left[y L_{2}^{(0)}(y)+\alpha(y)L_{1}^{(0)}(\alpha(y),y)\right]\right\}P_0^{'}(\zeta).
\end{split}
\end{gather}
Letting $y\rightarrow 1$ and then letting $x\rightarrow 1$ in \eqref{eq:3f9}, with L'H\^{o}pital's rule, we obtain
\begin{equation}\label{eq:3f10}
 \hspace*{\fill}Q_{1}^{(1)}(1,1,\zeta)=\frac{\rho_1}{1-\rho_1}[Q_{2}^{(1)}(1,\zeta)+Q_{3}^{(1)}(1,\zeta)]-\frac{\mu_3}{\mu_1}\frac{\rho_1(1-\rho_1-\rho_2)^2}{(1-\rho_1)^2}P_0^{'}(\zeta).\hspace*{\fill}
\end{equation}
Letting $y\rightarrow 1$ in \eqref{eq:3f11} and using L'H\^{o}pital's rule leads to
\begin{gather}\label{eq:3f12}
\begin{split}
&\frac{1}{1-\rho_1}\left[(1-\rho_1-\rho_2)Q_{2}^{(1)}(1,\zeta)-\rho_2Q_{3}^{(1)}(1,\zeta)\right]\\
 &=\left\{-\frac{\mu_3}{\mu_1}\left[\frac{\rho_1\rho_2}{1-\rho_1}+\frac{\rho_1\rho_2(1-\rho_1-\rho_2)}{(1-\rho_1)^2}\right]-\frac{\mu_3}{\mu_2}\rho_2\right\}P_0^{'}(\zeta).
\end{split}
\end{gather}
From \eqref{eq:3f10} and \eqref{eq:3f12}, we have
\begin{gather*}
\begin{split}\label{eq:3f13}
&(1-\rho_1-\rho_2)\left[Q_{1}^{(1)}(1,1,\zeta)+Q_{2}^{(1)}(1,\zeta)+Q_{3}^{(1)}(1,\zeta)\right]-Q_{3}^{(1)}(1,\zeta)\\
&=\frac{1}{1-\rho_1}\left[(1-\rho_1-\rho_2)Q_{2}^{(1)}(1,\zeta)-\rho_2Q_{3}^{(1)}(1,\zeta)\right]-\frac{\mu_3}{\mu_1}\frac{\rho_1(1-\rho_1-\rho_2)^2}{(1-\rho_1)^2}P_0^{'}(\zeta)\\
&=-\left[\frac{\mu_3}{\mu_1}\rho_1+\frac{\mu_3}{\mu_2}\rho_2\right]P_0^{'}(\zeta).
\end{split}\Box
\end{gather*}
\end{pf}

\subsection{Equating the second-order terms}\label{subsec4.3}
In this subsection we consider the sum of all $\mathcal{O}(\varepsilon^2)$ terms in equations \eqref{eq:36}-\eqref{eq:38} to determine $P_0(\zeta)$.

Taking the summation over all $x_1$ and $x_2$ of \eqref{eq:36}-\eqref{eq:38}, we get
\begin{gather}
\begin{split}\label{eq:3s1}
     & \mu_1\sum_{x_2=0}^{\infty}\phi_{1,(1,x_2)}(\zeta,\varepsilon)\\
      &=\lambda_1\sum_{x_2=1}^{\infty} \phi_{2,(0,x_2)}(\zeta,\varepsilon)+\lambda_1 \sum_{x_2=0}^{N-1}\phi_{3,(0,x_2)}(\zeta,\varepsilon)-\mu_3(1-\rho_1-\rho_2)\varepsilon\frac{\partial Q_1(1,1,\zeta,\varepsilon)}{\partial\zeta}\\
      &\ \ \ +\bigg[\omega\frac{\partial Q_1(1,1,\zeta,\varepsilon)}{\partial\zeta}+
      \frac{\mu_3(1-\rho_1-\rho_2)}{2}\frac{\partial^2Q_1(1,1,\zeta,\varepsilon)}{\partial\zeta^2}\bigg]\varepsilon^2+\mathcal{O}(\varepsilon^3),
\end{split}\\
\begin{split}\label{eq:3s2}
    &\lambda_1\sum_{x_2=1}^{\infty} \phi_{2,(0,x_2)}(\zeta,\varepsilon)+\mu_2\phi_{2,(0,1)}(\zeta,\varepsilon)\\
    &=\mu_1\sum_{x_2=1}^{\infty}\phi_{1,(1,x_2)}(\zeta,\varepsilon)+\lambda_2 \phi_{3,(0,N-1)}(\zeta,\varepsilon)-\mu_3(1-\rho_1-\rho_2)\varepsilon\frac{\partial Q_2(1,\zeta,\varepsilon)}{\partial\zeta}\\
    &\ \ \ +\left[\omega\frac{\partial Q_2(1,\zeta,\varepsilon)}{\partial\zeta}+
      \frac{\mu_3(1-\rho_1-\rho_2)}{2}\frac{\partial^2Q_2(1,\zeta,\varepsilon)}{\partial\zeta^2}\right]\varepsilon^2+\mathcal{O}(\varepsilon^3),
\end{split}\\
\begin{split}\label{eq:3s3}
      &\lambda_1 \sum_{x_2=0}^{N-1}\phi_{3,(0,x_2)}(\zeta,\varepsilon)+\lambda_2 \phi_{3,(0,N-1)}(\zeta,\varepsilon)\\
      &=\mu_1\phi_{1,(1,0)}(\zeta,\varepsilon)+\mu_2\phi_{2,(0,1)}(\zeta,\varepsilon)-\mu_3(1-\rho_1-\rho_2)\frac{\partial Q_3(1,\zeta,\varepsilon)}{\partial\zeta}+\mu_3\frac{\partial Q_3(1,\zeta,\varepsilon)}{\partial\zeta}\\
      &\ \ \ +\left[\omega\frac{\partial Q_3(1,\zeta,\varepsilon)}{\partial\zeta}+
      \frac{\mu_3(2-\rho_1-\rho_2)}{2}\frac{\partial^2Q_3(1,\zeta,\varepsilon)}{\partial\zeta^2}\right]\varepsilon^2 +\mathcal{O}(\varepsilon^3);
\end{split}
\end{gather}
Summing over \eqref{eq:3s1}-\eqref{eq:3s3}, we obtain
\begin{gather}
\begin{split}\label{eq:3s4}
      0=&\Bigg[-\mu_3(1-\rho_1-\rho_2)\left(\frac{\partial Q_1(1,1,\zeta,\varepsilon)}{\partial\zeta}+\frac{\partial Q_2(1,\zeta,\varepsilon)}{\partial\zeta}+\frac{\partial Q_3(1,\zeta,\varepsilon)}{\partial\zeta}\right)\\
      &+\mu_3\frac{\partial Q_3(1,\zeta,\varepsilon)}{\partial\zeta}\Bigg]\varepsilon+\Bigg[
      \frac{\mu_3(1-\rho_1-\rho_2)}{2
      }\Big(\frac{\partial^2Q_1(1,1,\zeta,\varepsilon)}{\partial\zeta^2}\\
      &+\frac{\partial^2Q_2(1,\zeta,\varepsilon)}{\partial\zeta^2}+\frac{\partial^2Q_3(1,\zeta,\varepsilon)}{\partial\zeta^2}\Big)+\omega \Big(\frac{\partial Q_1(1,1,\zeta,\varepsilon)}{\partial\zeta}\\
      &+\frac{\partial Q_2(1,\zeta,\varepsilon)}{\partial\zeta}+\frac{\partial Q_3(1,\zeta,\varepsilon)}{\partial\zeta}\Big)+\frac{\mu_3}{2}\frac{\partial^2Q_3(1,\zeta,\varepsilon)}{\partial\zeta^2}\Bigg]\varepsilon^2+\mathcal{O}(\varepsilon^3),
\end{split}
\end{gather}
Now taking the Taylor expansion \eqref{eq:380} of equation \eqref{eq:3s4}, we obtain
\begin{gather}
\begin{split}\label{eq:3s5}
      0=&\Big[\mu_3(1-\rho_1-\rho_2)P_0^{''}(\zeta)+\omega P_0^{'}(\zeta)+\mu_3Q_3^{'(1)}(1,\zeta)-\mu_3(1-\rho_1-\rho_2)\times\\
      &\left( Q_1^{'(1)}(1,1,\zeta)+Q_2^{'(1)}(1,\zeta)+Q_3^{'(1)}(1,\zeta)\right)\Big]\varepsilon^2+\mathcal{O}(\varepsilon^3)\\
      =&\left[\mu_3\left(1-\rho_1-\rho_2+\frac{\mu_3}{\mu_1}\rho_1+\frac{\mu_3}{\mu_2}\rho_2\right)P_0^{''}(\zeta)+\omega P_0^{'}(\zeta)\right]\varepsilon^2+\mathcal{O}(\varepsilon^3).
\end{split}
\end{gather}
In \eqref{eq:3s5}, the first equation follows from \eqref{eq:3181}-\eqref{eq:3201} and the second equation follows from Proposition \ref{pro1}.

From the above derivation procedure, we can conclude the following Proposition.
\begin{pro}\label{pro2}
After taking the summation over all $x_1$ and $x_2$ of the Taylor series of all perturbed balance equations \eqref{eq:36}-\eqref{eq:38}, the $\mathcal{O}(1)$ and $\mathcal{O}(\varepsilon)$ terms cancel and, moreover, equating the $\mathcal{O}(\varepsilon^2)$ terms yields the following differential equation for $P_0(\zeta)$:
\begin{equation*}\label{eq:3s70}
 \omega P_0^{'}(\zeta)=-\left[(1-\rho_1-\rho_2)+\frac{\mu_3}{\mu_1}\rho_1+\frac{\mu_3}{\mu_2}\rho_2\right]\mu_3P_0^{''}(\zeta).
\end{equation*}
\end{pro}

\subsection{The scaled number of customers in the critically loaded queue}\label{subsec4.4}
Now we can finally present the density of the scaled number of customers in $Q_3$, i.e. $P_0(\zeta)$. It can be obtained by combining the differential equation in Proposition \ref{pro2} with $\int_0^{\infty}P_0(\zeta)d\zeta=1$ that
\begin{equation*}\label{eq:3c1}
 P_0(\zeta)=\eta \mathrm{e}^{-\eta\zeta},
\end{equation*}
with $\frac{\omega}{\eta}=\left[1-\rho_1-\rho_2+\frac{\mu_3}{\mu_1}\rho_1+\frac{\mu_3}{\mu_2}\rho_2\right]\mu_3$.

As a special case, we may take $\omega=\mu_3$, which gives $\zeta=(1-\rho)X_3$, then
\begin{equation*}\label{eq:3c2}
 \frac{1}{\eta}=1-\rho_1-\rho_2+\frac{\mu_3}{\mu_1}\rho_1+\frac{\mu_3}{\mu_2}\rho_2.
\end{equation*}
By applying the multiclass distributional law of Bertsimas and Mourtzinou \cite{Bertsimas1997} it directly follows that the scaled waiting time at $Q_3$ follows an exponential distribution with parameter $\mu_3\eta$.

\subsection{Main result}\label{subsec4.5}
\begin{thm}\label{thm1}
For $\lambda_3=\mu_3(1-\rho_1-\rho_2)-\varepsilon\omega$, we have
\begin{equation*}\label{eq:3m}
\lim_{\varepsilon\downarrow 0}\mathbb{P}\{X_1\leq x_1,X_2\leq x_2,\varepsilon X_3\leq \zeta\}=\mathcal{L}(x_1,x_2)(1-\mathrm{e}^{-\eta\zeta}),
\end{equation*}
where $\mathcal{L}(\cdot,\cdot)$ is the joint cumulative distribution function(cdf) of the queue lengths of a preemptive priority polling system with N-policy vacation described in subsection \ref{subsec4.1}.
\end{thm}

The main result stated in Theorem \ref{thm1} can be interpreted as follows: in the heavy-traffic regime,
\begin{enumerate}[R1.]
  \item The queue lengths in the stable queues have the same distribution as that of a preemptive priority polling system with $N$-policy vacation.
  \item The scaled number of customers in the critically loaded queue is exponentially distributed with parameter $\eta$.
  \item The queue lengths in the stable queues and the (scaled) number of customers in the critically loaded queue are independent.
\end{enumerate}

For R1, since $Q_3$ is critically loaded, $Q_3$ would be visited during each cycle.
From the perspective of $Q_1$ and $Q_2$, the server goes on a vacation once the server goes to $Q_3$ when $Q_1$ and $Q_2$ are empty, and goes back to work once a type 1 customer arrives or there are $N$ type 2 customers queueing, which actually is an $N$-policy vacation.

For R2, we note that the total workload in the system equals the amount of workload in an M/G/1 queue with arrival rate $\lambda_1+\lambda_2+\lambda_3$ and hyperexponentially distributed service times, i.e. the service time is exponentially distributed with parameter $\mu_i$ with probability $\frac{\lambda_i}{\lambda_1+\lambda_2+\lambda_3}$, $i=1,2,3$.
Based on the heavy-traffic results for the M/G/1 queue (see \cite{Bertsimas1997}), the distribution of the scaled total workload converges to an exponential distribution with mean $\rho\mathbb{E}[R]$, where $R$ is a residual service time and
\begin{equation*}
 \hspace*{\fill}\mathbb{E}[R]=\frac{\frac{1}{\mu_1}\rho_1+\frac{1}{\mu_2}\rho_2+\frac{1}{\mu_3}\rho_3}{\rho}.\hspace*{\fill}
\end{equation*}
 In the heavy traffic, since almost all customers are located in $Q_3$, the total number of customers at this queue is also exponentially distributed with mean $\mu_3\left(\frac{1}{\mu_1}\rho_1+\frac{1}{\mu_2}\rho_2+\frac{1}{\mu_3}\rho_3\right)$. Since $\lambda_3\uparrow\mu_3(1-\rho_1-\rho_2)$, the scaled number of customers in $Q_3$ is exponentially distributed with parameter $\eta$.

Finally, R3 follows from the time-scale separation in the heavy traffic which implies that the dynamics of the stable queues evolve at a much faster time scale than the dynamics of the critically loaded queue. Since the amount of ``memory" of the stable queues asymptotically vanish compared to that of the critically loaded queue, the queue lengths in the stable queues are independent of the (scaled) number of customers in the critically loaded queue in the limit.

\begin{rmk}
From the above procedure, it is easy to see that, when there is a single critically loaded queue in the heavy traffic, the stable queues with threshold policies can always be transferred into a priority polling system with $N$-policy vacation.
\end{rmk}

\section{Exact tail asymptotics in Model \uppercase\expandafter{\romannumeral2} }\label{sec5}

In Section \ref{sec4}, we have derived the PGFs of the queue-length distributions of the stable queues, which have the same distributions as Model \uppercase\expandafter{\romannumeral2}.
As known, no closed-form expressions for the steady-state queue-length probabilities can be obtained.
In this section, we carry out a detailed analysis on the exact tail asymptotics for the stationary distributions in Model \uppercase\expandafter{\romannumeral2}, which provides us an approximation of the stable queues.

\subsection{Preliminary}\label{subsec5.1}
First we introduce some necessary notations. The marginal distributions for the high- and low-priority customers are denoted by $\pi_i^{(h)}$ and $\pi_j^{(l)}$, respectively. When $j>0$, we write $\pi_j^{(l)}=\pi_{1,j}^{(l)}+\pi_{2,j}^{(l)}$, where $\pi_{s,j}^{(l)}$ is the marginal distribution of the low-priority customers when the server is visiting $Q_s,s=1,2$. We denote the distribution of the total number of customers by $\pi_n^{(T)}$. Let $\lambda=\lambda_1+\lambda_2$ and $\overline{\rho}_1=\lambda/\mu_1$. Without loss of generality, throughout this section we assume that $\lambda_1+\lambda_2+\mu_1+\mu_2=1$. To completely derive the exact tail asymptotics, we first introduce the following notations:
\begin{align*}
&b_1=\frac{\lambda_2}{\lambda_2+(\sqrt{\mu_1}-\sqrt{\lambda_1})^2},\ \qquad\qquad\quad b_2=\frac{\lambda_2}{\lambda_2+(\sqrt{\mu_1}+\sqrt{\lambda_1})^2},\\
&\Delta(y)=(\lambda+\mu_1-\lambda_2y)^2-4\lambda_1\mu_1=\lambda_2^2(1-b_1y)(1-b_2y)/b_1b_2,\\
&x_1(y)=\frac{(\lambda+\mu_1-\lambda_2y)-\sqrt{\Delta(y)}}{2\lambda_1}=\alpha(y),\\
&x_2(y)=\frac{(\lambda+\mu_1-\lambda_2y)+\sqrt{\Delta(y)}}{2\lambda_1},\\
&xK(x,y)=-\lambda_1x^2+(\lambda+\mu_1-\lambda_2y)x-\mu_1=-\lambda_1(x-x_1(y))(x-x_2(y)),\\
&c_0=\frac{(\lambda+\mu_1)-\sqrt{(\lambda+\mu_1)^2-4\lambda_1\mu_1}}{2\mu_1},\ c_1=\frac{\lambda_2c_0}{\sqrt{(\lambda+\mu_1)^2-4\lambda_1\mu_1}},\\
&x_1=x_1(0)=\frac{c_0}{\rho_1},\qquad\qquad\qquad\qquad\quad x_2=x_2(0)=\frac{1}{c_0},\\
&F(y)=\lambda_2y^2-(1-2\mu_1+\mu_2)y+2\mu_2,\\
&T^*(y)=F(y)+y\sqrt{\Delta(y)},\ \ \ \ \qquad\qquad T(y)=F(y)-y\sqrt{\Delta(y)},\\
&\eta_1=\frac{(1-2\mu_1)+\sqrt{(1-2\mu_1)^2+4(\mu_1-\mu_2)\lambda_2}}{2\mu_2},\\
&\eta_2=\frac{(1-2\mu_1)-\sqrt{(1-2\mu_1)^2+4(\mu_1-\mu_2)\lambda_2}}{2\mu_2},\\
&T(y)T^*(y)=4\mu_2^2(1-y)(1-\eta_1y)(1-\eta_2y),\\
&a=\frac{1-\rho_1-\rho_2}{2\mu_2}\frac{\eta_1}{\eta_1-\eta_2},\ \ \qquad\ \ \ b=\frac{1-\rho_1-\rho_2}{2\mu_2}\frac{\eta_2}{\eta_2-\eta_1},\\
&D=(\lambda+\mu_1-2\sqrt{\lambda_1\mu_1})(\mu_1-\mu_2-\sqrt{\lambda_1\mu_1})+\lambda_2\mu_2.
\end{align*}

\subsection{The PGFs of the stationary queue-length distribution}\label{subsec5.2}
Define the following PGFs of the stationary queue-length distributions:
\begin{align*}
&\psi_j^{(0)}(x)=\sum_{i=1}^{\infty}\pi_{1,(i,j)}^{(0)}x^{i-1}, \ \ \ j=0,1,2,\ldots,\\
&L^{(l)}(y)=\sum_{n=0}^{\infty}\pi_n^{(l)}y^n,\qquad\quad\quad L^{(T)}(y)=\sum_{n=0}^{\infty}\pi_n^{(T)}y^n.
\end{align*}

Now we present some Propositions to give the exact expressions of the PGFs defined above.
\begin{pro}\label{pro5}
\begin{align}
&L_1^{(0)}(x,1)=\frac{\rho_1(1-\rho_1)}{1-\rho_1x},\label{eq:43}\\
&L_1^{(0)}(1,y)=\frac{\mu_2-\lambda_2y}{\lambda_2}L_2^{(0)}(y)-L_3^{(0)}(y),\label{eq:44}\\
&L_1^{(0)}(y,y)=\frac{\lambda y-\mu_2}{\mu_1(1-\bar{\rho}_1y)}L_2^{(0)}(y)+\frac{\lambda}{\mu_1(1-\bar{\rho}_1y)}L_3^{(0)}(y),\label{eq:45}\\
&L_1^{(0)}(x,0)=\frac{c_0}{1-c_0x}L_3^{(0)}(0),\label{eq:451}
\end{align}
where $L_1^{(0)}(x,y)$ and $L_2^{(0)}(y)$ are expressed in \eqref{eq:318} and \eqref{eq:319} respectively.
\end{pro}

\begin{pf}
Adding \eqref{eq:a6} to \eqref{eq:a5}, we have
\begin{equation}\label{eq:431}
xK(x,y)L_1^{(0)}(x,y) =y[\lambda_1 x-a(y)]L_2^{(0)}(y)+[\lambda_1 (x-1)+\lambda_2(y-1)]L_3^{(0)}(y).
\end{equation}
Then, letting $y\rightarrow 1$, $x\rightarrow 1$, $x\rightarrow y$ and $y\rightarrow 0$, respectively, we get equations \eqref{eq:43}-\eqref{eq:451}.\hspace*{\fill}$\Box$
\end{pf}

\begin{pro}\label{pro6}
\begin{align}
&\psi_0^{(0)}(x)=\frac{c_0}{1-c_0x}L_{3}^{(0)}(0),\label{eq:452}\\
&\psi_j^{(0)}(x)=\frac{a_j}{1-c_0x}+\frac{\lambda_2c_0x}{\lambda_1(1-c_0x)}\frac{\psi_{j-1}^{(0)}(x)-\psi_{j-1}^{(0)}(x_1)}{x-x_1}, \ j=1,2,\ldots\label{eq:453}
\end{align}
where $a_j=\frac{c_0}{\lambda_1}\left[\lambda_2\psi_{j-1}^{(0)}(x_1)+\lambda_1\left(\pi_{2,(0,j)}^{(0)}+\pi_{3,(0,j)}^{(0)}\delta(j<N)\right)\right]$.
\end{pro}

\begin{pf}
Equation \eqref{eq:452} is obvious since $\psi_0^{(0)}(x)=L_1^{(0)}(x,0)$. Using the PGFs to rewrite balance equation \eqref{eq:39}, we obtain
\begin{equation}\label{eq:4531}
\psi_j^{(0)}(x)=\frac{\lambda_2x\psi_{j-1}^{(0)}(x)+\lambda_1x\left(\pi_{2,(0,j)}^{(0)}+\pi_{3,(0,j)}^{(0)}\delta(j<N)\right)-\mu_1\psi_j^{(0)}(0)}{-\lambda_1(x-x_1)(x-x_2)}.
\end{equation}
Note that $x_1<1$ and $\psi_j^{(0)}(x)$ is analytic inside the unit circle, which implies that $x_1$ is also a zero of the numerator of the righthand side of \eqref{eq:4531}. Therefore,
\begin{equation}\label{eq:4532}
\lambda_2x_1\psi_{j-1}^{(0)}(x_1)+\lambda_1x_1\left(\pi_{2,(0,j)}^{(0)}+\pi_{3,(0,j)}^{(0)}\delta(j<N)\right)=\mu_1\psi_j^{(0)}(0).
\end{equation}
Taking \eqref{eq:4532} into the numerator of the right hand side of \eqref{eq:4531} yields
\begin{gather*}\label{eq:4533}
\begin{split}
\psi_j^{(0)}(x)=&\frac{\left[\lambda_2\psi_{j-1}^{(0)}(x_1)+\lambda_1\left(\pi_{2,(0,j)}^{(0)}+\pi_{3,(0,j)}^{(0)}\delta(j<N)\right)\right](x-x_1)}{-\lambda_1(x-x_1)(x-x_2)}\\
&+\frac{\lambda_2x\left(\psi_{j-1}^{(0)}(x)-\psi_{j-1}^{(0)}(x_1)\right)}{-\lambda_1(x-x_1)(x-x_2)}.
\end{split}
\end{gather*}
Since $x_2=\frac{1}{c_0}$, \eqref{eq:453} can be obtained by simplifying the above equation.\hspace*{\fill}$\Box$
\end{pf}

\begin{pro}\label{pro7}
\begin{equation*}\label{eq:46}
L_2^{(0)}(y)=\left[\frac{aT^*(y)}{1-\eta_1y}+\frac{bT^*(y)}{1-\eta_2y}\right]\iota(y)\beta(y),
\end{equation*}
with $\iota(y)=\frac{\mu_1-\lambda+\lambda_2y-\sqrt{\Delta(y)}}{2\mu_2(y-1)}$.
\end{pro}

\begin{pf}
Simplifying \eqref{eq:319}, we get
\begin{gather*}
\begin{split}
L_2^{(0)}(y)&=\frac{\lambda_1(x_1(y)-1)+\lambda_2(y-1)}{ya(y)-\lambda_1yx_1(y)}L_3^{(0)}(y)\\
&=(1-\rho_1-\rho_2)\frac{2T^*(y)\mu_2(1-y)}{T(y)T^*(y)}\frac{\lambda_1(x_1(y)-1)+\lambda_2(y-1)}{\mu_2(y-1)}\beta(y)\\
&=\frac{1-\rho_1-\rho_2}{2\mu_2}\frac{T^*(y)}{(1-\eta_1y)(1-\eta_2y)}\iota(y)\beta(y)\\
&=\left[\frac{aT^*(y)}{1-\eta_1y}+\frac{bT^*(y)}{1-\eta_2y}\right]\iota(y)\beta(y).
\end{split}\Box
\end{gather*}
\end{pf}

\begin{pro}\label{pro8}
\begin{equation*}\label{eq:461}
L^{(T)}(y)=\left[a\frac{T^*(y)}{1-\eta_1y}+b\frac{T^*(y)}{1-\eta_2y}\right]\kappa(y)\beta(y),
\end{equation*}
with $\kappa(y)=\frac{2\mu_2(1-y)-(1-\mu_1y)T(y)}{2\mu_1\mu_2y(1-y)(1-\bar{\rho}_1y)}$.
\end{pro}
\begin{pf}
By the definition of $L^{(T)}(y)$, we have
\begin{gather*}
\begin{split}
L^{(T)}(y)&=yL_1^{(0)}(y,y)+yL_2^{(0)}(y)+L_3^{(0)}(y)\\
&=\frac{1}{1-\bar{\rho}_1y}L_3^{(0)}(y)+\frac{(\mu_1-\mu_2)y}{\mu_1(1-\bar{\rho}_1y)}L_2^{(0)}(y)\\
&=\left[a\frac{T^*(y)}{1-\eta_1y}+b\frac{T^*(y)}{1-\eta_2y}\right]\kappa(y)\beta(y),
\end{split}
\end{gather*}
where the second equation follows from the expression \eqref{eq:45} and the last follows the same idea used in Proposition \ref{pro7}.\hspace*{\fill}$\Box$
\end{pf}

\subsection{Analysis of singularities and asymptotic expansions}\label{subsec5.2}
Along the same idea used for the classical priority model in \cite{Li2009}, asymptotics of the coeffients  are obtained using the following Tauberian-like theorem, which is Corollary 2 given in \cite{Flajolet1990}. For a function $f(y)$ that is analytic at $y=0$, we denote the coefficient of $y^k$ in the Taylor expression of $f(y)$ by $C_k[f(y)]$.

For the compactness, we omit all the proofs in this subsection, which can be referred to \cite{Li2009}.

\begin{lem}[Flayolet and Odlyzko]\label{lem3}
Assume that $f(z)$ is analytic in $\Delta(\phi,\varepsilon)=\{z:|z|\leq 1+\varepsilon, |Arg(z-1)|\geq \phi \ \text{for}\  \varepsilon>0 \ \text{and}\  0<\phi<\pi/2\}$ except at $z=1$ and
\begin{equation*}
f(z)\sim K(1-z)^s \qquad \text{as}\  z\rightarrow 1 \qquad \text{in}\  \Delta(\phi,\varepsilon).
\end{equation*}
Then as $n\rightarrow\infty$:
\begin{enumerate}[1.]
  \item If $s \notin \{0,1,2,\ldots\}$,
  \begin{equation*}
f_n\sim \frac{K}{\Gamma(-s)}n^{-s-1}.
\end{equation*}
  \item If $s$ is a nonnegative integer, then
  \begin{equation*}
f_n=o(n^{-s-1}).
\end{equation*}
\end{enumerate}
\end{lem}

The key goal is to locate the dominant singularity, which determines the decay and to characterize the nature of the dominant singularity, which determines the prefactor and the singularity coefficient.

Define
\begin{multline*}
\tilde{\Delta}(\phi,\varepsilon,a)=\{z:|az|\leq 1+\varepsilon, |Arg(az-1)|\geq \phi \ \text{for} \  0<a<1,\\
\  \varepsilon>0 \ \text{and}\  0<\phi<\pi/2\}-\{1/a\}.
\end{multline*}

\begin{lem}
For the non-unit zeros $1/\eta_1$ and $1/\eta_2$, we have
\begin{enumerate}[1.]
  \item Both $1/\eta_1$ and $1/\eta_2$ are real.
  \item $\eta_1>0$.
  \item $\eta_1>\eta_2$, and $\eta_1<|\eta_2|$ implies $\eta_2<0$.
  \item $\eta_2<0$, $\eta_2=0$ or $\eta_2>0$ if and only if $\mu_2<\mu_1$, $\mu_2=\mu_1$ or $\mu_2>\mu_1$, respectively.
  \item $\eta_2\neq b_1$, and either $T*(1/\eta_2)=0$ or $|\eta_2|< b_1$.
\end{enumerate}
\end{lem}

\begin{lem}[Key Lemma]
There are three cases for the dominant singularity of $L_2^{(0)}(y)$:
\begin{enumerate}[1.]
  \item If $D>0$, then $1<1/\eta_1<1/b_1$ and $1/\eta_1$ is a zero of $T(y)$ (but not $T^*(y)$), and therefore $1/\eta_1$ is the dominant singularity of $L_2^{(0)}(y)$, which is a simple pole.
  \item If $D=0$, then $1<1/\eta_1=1/b_1$ and $1/\eta_1$ is a zero of $T(y)$ and $T^*(y)$, and therefore $1/b_1$ is the dominant singularity of $L_2^{(0)}(y)$, which is both a branch point and a simple pole.
  \item If $D<0$, then $1<1/\eta_1<1/b_1$ and $1/\eta_1$ is a zero of $T^*(y)$ (but not $T(y)$), and therefore $1/b_1$ is the dominant singularity of $L_2^{(0)}(y)$, which is a branch point.
\end{enumerate}
\end{lem}

\begin{pro}\label{pro9}
If $\eta$ satisfies: (\romannumeral 1) $\eta\neq0$; (\romannumeral 2) $\eta\neq b_1$; (\romannumeral 3) $|\eta|<b_1$ or $T^*(\eta)=0$, then for $\eta=\eta_i$, $i=1,2$,
\begin{equation*}\label{eq:4e1b}
C_n\left[\frac{T^*(y)}{1-\eta y}\iota(y)\beta(y)\right]\sim b_1\beta(1/b_1)\sigma(\eta)n^{-3/2}b_1^n,
\end{equation*}
with $\sigma(\eta)=\frac{K(\eta)}{b_1\sqrt{\pi}}$ and $K(\eta)=\frac{\lambda_2b_1\sqrt{1-b_2/b_1}}{2\sqrt{b_1b_2}(\eta-b_1)}$.
\end{pro}

 \begin{pro}\label{pro10}
If $\bar{\rho}_1\geq 1$ and $\eta$ satisfy: (\romannumeral 1)
$\eta\neq0$; (\romannumeral 2) $\eta\neq b_1$; (\romannumeral 3)
$|\eta|<b_1$ or $T^*(\eta)=0$, then for $\eta=\eta_i$, $i=1,2$,
\begin{equation*}\label{eq:4e2b}
C_n\left[\frac{T^*(y)}{1-\eta y}\kappa(y)\beta(y)\right]\sim \beta(1/b_1)\sigma_1(\eta)n^{-3/2}b_1^n,
\end{equation*}
with $\sigma_1(\eta)=\frac{K_1(\eta)}{b_1\sqrt{\pi}}$ and $K_1(\eta)=\frac{\lambda_2b_1\sqrt{1-b_2/b_1}\big[(1-\mu_1/b_1)\big((F(1/b_1)+1\big)-2\mu_2(1-1/b_1)\big]}{4\mu_1\mu_2\sqrt{b_1b_2}(1-\eta/b_1)(1-1/b_1)(1-\bar{\rho_1}/b_1)}$.
\end{pro}

\subsection{Main results of exact tail asymptotics}\label{subsec5.3}
In this subsection, we provide a complete exact tail asymptotics of the stationary distributions(the joint and marginal queue lengths and the total number of customers ) by using the Tauberian-like Theorem to the related generating functions.
\begin{thm}
The exact tail asymptotics in the marginal stationary distribution $\pi_n^{(h)}$ of the high-priority queue is given by
\begin{equation*}
\pi_n^{(h)}\sim (1-\rho_1)\rho_1^n.
\end{equation*}
The decay rate in the marginal distribution for the high-priority queue is $\rho_1$.
\end{thm}

\begin{pf}
It is a direct consequence of the Taylor expansion of \eqref{eq:43}.\hspace*{\fill}$\Box$
\end{pf}

\begin{thm}
The exact tail asymptotics in the joint stationary distribution along the high-priority queue is characterized by: for a fixed number $j\geq 0$ of low-priority customers,
\begin{equation*}
\pi_{1,(n,j)}^{(0)}\sim \beta(0)(1-\rho_1-\rho_2)\left(\frac{c_1^j}{j!}\right)n^jc_0^{n-j}.
\end{equation*}
\end{thm}

\begin{pf}
First, by the induction, we prove
\begin{equation}\label{47}
\psi_j^{(0)}(x)\sim c_0 \beta(0)(1-\rho_1-\rho_2)\left(\frac{c_1}{c_0}\right)^j\frac{1}{(1-c_0x)^{j+1}},\ j\geq 0, \ \text{as}\  c_0x\rightarrow 1.
\end{equation}
 It is true for $j=0$ since $\psi_0^{(0)}(x)=\frac{c_0}{1-c_0x}L_{3}^{(0)}(0)$. Assume that \eqref{47} is true for $j=k$, we then show it is true for $j=k+1$. Rewrite equation \eqref{eq:453} as
\begin{equation*}
\psi_{k+1}^{(0)}(x)=\frac{a_{k+1}}{1-c_0x}+\frac{\lambda_2c_0x}{\lambda_1(1-c_0x)}\frac{\psi_{k}^{(0)}(x)-\psi_{k}^{(0)}(x_1)}{x-x_1},
\end{equation*}
where $a_{k+1}$ is a constant. Note that $\frac{\lambda_2}{\lambda_1}\frac{c_0}{1-x_1c_0}=\frac{c_0}{c_1}$. Hence,
\begin{equation*}
\lim_{c_0x\rightarrow 1}\frac{\psi_{k+1}^{(0)}(x)}{(1-c_0x)^{-(k+2)}}=c_0 \beta(0)(1-\rho_1-\rho_2)\left(\frac{c_1}{c_0}\right)^{k+1},
\end{equation*}
which is equivalent to \eqref{47}. Therefore, \eqref{47} is true for all $j\geq 0$.

Applying Lemma \ref{lem3} to \eqref{47}, we have
\begin{equation*}
\frac{C_n[\psi_j^{(0)}(x)]}{c_0^n}\sim c_0 L_{3}^{(0)}(0)\left(\frac{c_1}{c_0}\right)^j\frac{n^{(j+1)-1}}{\Gamma(j+1)}=c_0 L_{3}^{(0)}(0)\left(\frac{c_1}{c_0}\right)^j\frac{n^j}{j!}, \ j\geq 0,
\end{equation*}
that is
\begin{equation*}
\pi_{1,(n+1,j)}^{(0)}\sim L_{3}^{(0)}(0)\left(\frac{c_1^j}{j!}\right)n^jc_0^{n+1-j}, \ j\geq 0,
\end{equation*}
which completes the proof.\hspace*{\fill}$\Box$
\end{pf}

\begin{thm}\label{thm4}
The exact tail asymptotics in the joint stationary distribution along the low-priority queue is characterized by: for a fixed number $i\geq 0$ of high-priority customers,
  \begin{enumerate}[1.]
    \item (Exact geometric decay) In the region of $D>0$,
    \begin{equation*}\label{48}
\pi_{2,(i,n)}^{(0)}\sim C_{2,l,1}[u(\eta_1)]^i\eta_1^n.
\end{equation*}
    \item (Geometric decay with prefactor $n^{-1/2}$) In the region of $D=0$,
    \begin{equation*}\label{49}
\pi_{2,(i,n)}^{(0)}\sim C_{2,l,2}(\sqrt{\rho_1})^in^{-1/2}b_1^n.
\end{equation*}
    \item (Geometric decay with prefactor $n^{-3/2}$) In the region of $D<0$,
    \begin{equation*}\label{410}
\pi_{2,(i,n)}^{(0)}\sim C_{2,l,2}(1+i\widetilde{B})(\sqrt{\rho_1})^in^{-3/2}b_1^n.
\end{equation*}
  \end{enumerate}
Here $C_{2,l,1}$, $C_{2,l,2}$, $C_{2,l,3}$, $u(\eta)$ and $\widetilde{B}$ are given below:
\begin{align*}
&C_{2,l,1}=2aF(1/\eta_1)\beta(1/\eta_1),\\
&C_{2,l,2}=\frac{a\lambda_2\sqrt{1-b_2/b_1}}{\sqrt{\pi}b_1\sqrt{b_1b_2}}\beta(1/b_1),\\
&C_{2,l,3}=[a\sigma(\eta_1)+b\sigma(\eta_2)]\beta(1/b_1),\\
&u(\eta)=\frac{1-\mu_2-(\lambda_2/\eta)-\sqrt{\left[1-\mu_2-(\lambda_2/\eta)\right]^2-4\lambda_1\mu_1}}{2\mu_1},\\
&\widetilde{B}=\frac{\mu_2-\mu_1-\mu_2b_1+\sqrt{\lambda_1\mu_1}}{\sqrt{\lambda_1\mu_1}}.
\end{align*}
\end{thm}

\begin{pf}
In the case of $i=0$,
\begin{enumerate}
  \item if $D>0$, then $T(\frac{1}{\eta_1})=0$, and we can prove $\iota(\frac{1}{\eta_1})=\eta_1$. Hence,
  \begin{gather*}
  \begin{split}
&\lim_{\eta_1y\rightarrow 1}\left[\frac{L_2^{(0)}(y)}{(1-\eta_1y)^{-1/2}}\right]\\
&=a\lim_{\eta_1y\rightarrow 1}T^*(y)\iota(y)\beta(y)+b\lim_{\eta_1y\rightarrow 1}\left[\frac{(1-\eta_1y)T^*(y)}{1-\eta_2y}\iota(y)\beta(y)\right]\\
&=aT^*(1/\eta_1)\iota(1/\eta_1)\beta(1/\eta_1)=2aF(1/\eta_1)\eta_1\beta(1/\eta_1).
\end{split}
\end{gather*}
Clearly, $L_2^{(0)}(y)$ is analytic in $\tilde{\Delta}(\phi,\varepsilon,\eta_1)$. By Lemma \ref{lem3}, we obtain
\begin{equation*}
\pi_{2,(0,n+1)}^{(0)}\sim C_{2,l,1}\eta_1^{n+1}.
\end{equation*}

  \item if $D=0$, then $T(\frac{1}{b_1})=T(\frac{1}{\eta_1})=0$, hence, $\iota(\frac{1}{\eta_1})=b_1$ and
\begin{gather*}
\begin{split}
&\frac{T^*(y)}{1-\eta_1y}\iota(y)\beta(y)\\
&= \frac{F(y)-F(1/b_1)}{1-b_1y}\iota(y)\beta(y)+\frac{y\sqrt{\Delta(y)}}{1-b_1y}\iota(y)\beta(y)\\
&\sim \frac{\rho_2F'(1/b_1)\sqrt{1-b_2/b_1}}{2(1-b_1)\sqrt{b_1b_2}} \sqrt{1-b_1y}+\frac{\lambda_2\sqrt{1-b_2/b_1}\beta(1/b_1)}{\sqrt{b_1b_2}\sqrt{1-b_1y}}\\
&\sim \frac{\lambda_2\sqrt{1-b_2/b_1}\beta(1/b_1)}{\sqrt{b_1b_2}\sqrt{1-b_1y}}.
\end{split}
\end{gather*}

Since $\frac{T^*(y)}{1-\eta_1y}\iota(y)\beta(y)$ is analytic in $\tilde{\Delta}(\phi,\varepsilon,b_1)$, applying Lemma \ref{lem3}, we get
\begin{equation*}
C_n\left[\frac{T^*(y)}{1-\eta_1y}\iota(y)\beta(y)\right]\sim \frac{\lambda_2\sqrt{1-b_2/b_1}\beta(1/b_1)}{\sqrt{b_1b_2}\sqrt{\pi}}n^{-3/2}b_1^{n+1}.
\end{equation*}
While with Proposition \ref{pro9}, we have
\begin{equation*}
C_n\left[\frac{T^*(y)}{1-\eta_2y}\iota(y)\beta(y)\right]\sim \beta(1/b_1)\sigma(\eta_2)n^{-1/2}b_1^{n+1}.
\end{equation*}
Combining the above two asymptotics gives
\begin{equation*}
\pi_{2,(0,n+1)}^{(0)}\sim C_{2,l,2}n^{-1/2}b_1^{n+1}.
\end{equation*}

\item if $D<0$, the conclusion is a direct consequence of Proposition \ref{pro9}.
\end{enumerate}

In the case of $i>0$, the theorem can be proved by induction on $i$.
\begin{enumerate}
  \item if $D>0$, for $i=1$, the balance equation is
\begin{equation*}
\mu_1 \frac{\pi_{1,(1,n)}^{(0)}}{\eta_1^n}=(\lambda_1+\lambda_2 +\mu_2)\frac{\pi_{2,(0,n)}^{(0)}}{\eta_1^n}-\frac{\lambda_2}{\eta_1}\frac{\pi_{2,(0,n-1)}^{(0)}}{\eta_1^{n-1}}
                           -\mu_2\eta_1\frac{\pi_{2,(0,n+1)}^{(0)}}{\eta_1^{n+1}} .
\end{equation*}
 It is easy to see that $u(\eta)$ is the root of the equation with smaller module: $\mu_1[t(\eta)]^2-[1-\mu_2-\lambda_2/\eta]t(\eta)+\lambda_1=0$. Since $T(\frac{1}{\eta})=0$, we have $u(\eta)=\frac{1-\mu_1-\mu_2\eta-\lambda_2/\eta}{\mu_1}$. Therefore, we obtain
\begin{equation*}
\pi_{1,(1,n)}^{(0)}\sim C_{2,l,1}A_1\eta_1^n,
\end{equation*}
where $A_1=u(\eta_1)$. Assume that for $i\leq k$,
\begin{equation*}
\pi_{1,(i,n)}^{(0)}\sim C_{2,l,1}A_i\eta_1^n.
\end{equation*}
 Based on the balance equation
\begin{gather*}
\mu_1 \pi_{1,(2,n)}^{(0)}=(\lambda_1+\lambda_2 +\mu_2)\pi_{1,(1,n)}^{(0)}-\lambda_2\pi_{1,(1,n-1)}^{(0)}
                           -\lambda_1\pi_{2,(0,n)}^{(0)},\\
\mu_1 \pi_{1,(k+1,n)}^{(0)}=(\lambda_1+\lambda_2 +\mu_2)\pi_{1,(k,n)}^{(0)}-\lambda_2\pi_{1,(k,n-1)}^{(0)}
                           -\lambda_1\pi_{1,(k-1,n)}^{(0)},
\end{gather*}
and the inductive assumption $\frac{\pi_{1,(k+1,n)}^{(0)}}{\eta_1^n}\rightarrow C_{2,l,1}A_{k+1}$, we have
\begin{equation*}
\mu_1 A_{k+1}=(\lambda_1+\lambda_2 +\mu_2-\frac{\lambda_2}{\eta_1})A_k-\lambda_1A_{k-1}, \ k=1,2,3,\ldots
\end{equation*}
with $A_0=1$ and $A_1=u(\eta_1)$. Solving this difference equation leads to
\begin{equation*}
A_{k}=[u(\eta_1)]^k, \ k=0,1,2,\ldots
\end{equation*}
which gives the conclusion.

  \item if $D=0$, the proof is similar to that for case 1.

\item if $D<0$, then $u(b_1)=\sqrt{\rho_1}$. Along the same idea in the proof of case 1, we get a difference equation
    \begin{equation*}
 A_{k+1}=2\sqrt{\rho_1}A_k-\rho_1A_{k-1}, \ k=1,2,3,\ldots
\end{equation*}
with $A_0=1$ and $A_1=u(b_1)$. Solving the equation yields the conclusion.
\end{enumerate}\hspace*{\fill}$\Box$
\end{pf}

\begin{thm}\label{thm5}
The exact tail asymptotics in the marginal stationary distribution $\pi_n^{(l)}$ of the low-priority queue is given by
\begin{equation*}
\pi_n^{(l)}=\frac{\mu_2}{\lambda_2}\pi_{2,(0,n+1)}.
\end{equation*}
\end{thm}
\begin{pf}
It is clear since $L^{(l)}(y)=L_1^{(0)}(1,y)+yL_2^{(0)}(y)=\frac{\mu_2}{\lambda_2}L_2^{(0)}(y)-L_3^{(0)}(y)$.\hspace*{\fill}$\Box$
\end{pf}

\begin{thm}\label{thm6}
The exact tail asymptotics in the stationary distribution $\pi_n^T$ of total number of customers in the system is characterized below:

If $\mu_1=\mu_2$, then
\begin{equation*}
\pi_n^T=\beta\left(\frac{1}{1-\rho_1-\rho_2}\right)(1-\rho_1-\rho_2)(\rho_1+\rho_2)^n,\quad n=0,1,2,\ldots.
\end{equation*}

If $\mu_1\neq\mu_2$, then
\begin{enumerate}
  \item In the region of $D>0$, three cases exist:
  \begin{enumerate}[a)]
    \item If (\romannumeral 1) $\bar{\rho}_1\geq 1$; or (\romannumeral 2) $\bar{\rho}_1< 1$ and $\bar{\rho}_1<\eta_1$, then
    \begin{equation*}
\pi_n^T\sim C_{t,1a}\eta_1^{n}.
\end{equation*}
    \item If $\bar{\rho}_1< 1$ and $\bar{\rho}_1>\eta_1$, then
    \begin{equation*}
\pi_n^T\sim C_{t,1b}(\bar{\rho}_1)^{n}.
\end{equation*}
    \item If $\bar{\rho}_1< 1$ and $\bar{\rho}_1=\eta_1$, then
    \begin{equation*}
\pi_n^T\sim C_{t,1c}n\eta_1^{n}.
\end{equation*}
  \end{enumerate}
  \item In the region of $D=0$, two cases exist:
  \begin{enumerate}[a)]
    \item If $\bar{\rho}_1\geq 1$, then
    \begin{equation*}
\pi_n^T\sim C_{t,2a}n^{-1/2}b_1^{n}.
\end{equation*}
    \item If $\bar{\rho}_1< 1$, then
    \begin{equation*}
\pi_n^T\sim C_{t,2b}(\bar{\rho}_1)^{n}.
\end{equation*}
  \end{enumerate}
  \item In the region of $D<0$, three cases exist:
  \begin{enumerate}[a)]
    \item If $\bar{\rho}_1\geq 1$, then
    \begin{equation*}
\pi_n^T\sim C_{t,3a}n^{-3/2}b_1^{n}.
\end{equation*}
    \item If $\bar{\rho}_1< 1$ and $\bar{\rho}_1\neq \sqrt{\rho_1}$,then
    \begin{equation*}
\pi_n^T\sim C_{t,3b}(\bar{\rho}_1)^{n}.
\end{equation*}
\item If $\bar{\rho}_1< 1$ and $\bar{\rho}_1= \sqrt{\rho_1}$,then $\bar{\rho}_1=b_1\neq \eta_1$ and
    \begin{equation*}
\pi_n^T\sim C_{t,3c}(\bar{\rho}_1)^{n}.
\end{equation*}
\end{enumerate}
\end{enumerate}
Here $C_{t,1a}$, $C_{t,1b}$, $C_{t,1c}$, $C_{t,2a}$, $C_{t,2b}$, $C_{t,3a}$, $C_{t,3b}$ and $C_{t,3c}$ are given below:
\begin{align*}
&C_{t,1a}=\frac{(\mu_1-\mu_2)\eta_1}{\mu_1(\eta_1-\bar{\rho}_1)}C_{2,l,1},\\
&C_{t,1b}=C_{t,2b}=C_{t,3b}=C_{t,3c}=\frac{(\mu_1-\mu_2)}{\mu_1}\frac{1}{\bar{\rho}_1}L_2^{(0)}(\frac{1}{\bar{\rho}_1})+L_3^{(0)}(\frac{1}{\bar{\rho}_1}),\\
&C_{t,1c}=\frac{(\mu_1-\mu_2)}{\mu_1}C_{2,l,1},\\
&C_{t,2a}=\frac{\kappa(1/b_1)}{b_1}C_{2,l,2},\\
&C_{t,3a}=[a\sigma_1(\eta_1)+b\sigma_1(\eta_2)]\beta(1/b_1).
\end{align*}
\end{thm}

\begin{pf}
If $\mu_1=\mu_2$, then $L^{(T)}(y)=\frac{1}{1-\bar{\rho}_1y}L_3^{(0)}(y)$ and $\bar{\rho}_1=\rho_1+\rho_2$. Hence, the conclusion is true. Now we consider the case $\mu_1\neq \mu_2$.
\begin{enumerate}
   \item In the region of $D>0$,
  \begin{enumerate}[a)]
    \item If (\romannumeral 1) $\bar{\rho}_1\geq 1$; or (\romannumeral 2) $\bar{\rho}_1< 1$ and $\bar{\rho}_1<\eta_1$, then
  \begin{gather*}
  \begin{split}
\lim_{\eta_1y\rightarrow 1}\left[\frac{L^{(T)}(y)}{(1-\eta_1y)^{-1}}\right]&=a\lim_{\eta_1y\rightarrow 1}\frac{(\mu_1-\mu_2)y}{\mu_1(1-\bar{\rho}_1y)}T^*(y)\iota(y)\beta(y)\\
&=\frac{(\mu_1-\mu_2)\eta_1}{\mu_1(\eta_1-\bar{\rho}_1)}C_{2,l,1}.
\end{split}
\end{gather*}
Since $L^{(T)}(y)$ is analytic in $\tilde{\Delta}(\phi,\varepsilon,\eta_1)$, applying Lemma \ref{lem3}, we get
\begin{equation*}
\pi_{n}^{(T)}\sim C_{t,1a}\eta_1^{n}.
\end{equation*}
\item If $\bar{\rho}_1< 1$ and $\bar{\rho}_1>\eta_1$, then $L^{(T)}(y)$ is analytic in $\tilde{\Delta}(\phi,\varepsilon,\bar{\rho}_1)$ and
\begin{equation*}
\lim_{\bar{\rho}_1y\rightarrow 1}\left[\frac{L^{(T)}(y)}{(1-\bar{\rho}_1y)^{-1}}\right]=C_{t,1b}.
\end{equation*}
By Lemma \ref{lem3}, we have
\begin{equation*}
\pi_{n}^{(T)}\sim C_{t,1b}(\bar{\rho}_1)^{n}.
\end{equation*}
\item If $\bar{\rho}_1< 1$ and $\bar{\rho}_1=\eta_1$, then $L^{(T)}(y)$ is analytic in $\tilde{\Delta}(\phi,\varepsilon,\eta_1)$ and
\begin{equation*}
\lim_{\eta_1y\rightarrow 1}\left[\frac{L^{(T)}(y)}{(1-\eta_1y)^{-2}}\right]=\frac{(\mu_1-\mu_2)}{\mu_1}C_{2,l,1}.
\end{equation*}
By Lemma \ref{lem3}, we obtain
\begin{equation*}
\pi_{n}^{(T)}\sim C_{t,1c}n\eta_1^{n}.
\end{equation*}
\end{enumerate}

\item In the region of $D=0$, two cases exist:
  \begin{enumerate}[a)]
    \item If $\bar{\rho}_1\geq 1$, by Proposition \ref{pro8},
    \begin{equation*}
L^{(T)}(y)=\left[a\frac{T^*(y)}{1-\eta_1y}+b\frac{T^*(y)}{1-\eta_2y}\right]\kappa(y)\beta(y).
\end{equation*}
    Similarly to the case of $D=0$ in Theorem \ref{thm5}, we have
    \begin{equation*}
\frac{T^*(y)}{1-\eta_1y}\kappa(y)\beta(y)\sim \frac{\lambda_2\sqrt{1-b_2/b_1}\kappa(1/b_1)\beta(1/b_1)}{\sqrt{b_1b_2}\sqrt{1-b_1y}}.
\end{equation*}
In addition, $\frac{T^*(y)}{1-\eta_1y}\kappa(y)\beta(y)$ is analytic in $\tilde{\Delta}(\phi,\varepsilon,b_1)$. Hence,
\begin{equation*}
C_n\left[\frac{T^*(y)}{1-\eta_1y}\kappa(y)\beta(y)\right]\sim \frac{\lambda_2\sqrt{1-b_2/b_1}\kappa(1/b_1)\beta(1/b_1)}{\sqrt{b_1b_2}\sqrt{\pi}}n^{-1/2}b_1^{n}.
\end{equation*}
While with Proposition \ref{pro10}, we have
\begin{equation*}
C_n\left[\frac{T^*(y)}{1-\eta_2y}\kappa(y)\beta(y)\right]\sim \beta(1/b_1)\sigma_1(\eta_2)n^{-3/2}b_1^{n+1}.
\end{equation*}

Combining the above two asymptotics leads to
\begin{equation*}
\pi_n^T\sim C_{t,2a}n^{-1/2}b_1^{n}.
\end{equation*}
    \item If $\bar{\rho}_1< 1$, then $\bar{\rho}_1>b_1$. This can be proved by contradiction: if $\bar{\rho}_1=b_1$, then $\bar{\rho}_1=\sqrt{\rho_1}$, which follows from $\bar{\rho}_1-b_1=\frac{-(\bar{\rho}_1-\sqrt{\rho_1})^2}{\bar{\rho}_1+1-2\sqrt{\rho_1}}$. After some manipulations, we get $D=\mu_1(1-\sqrt{\rho_1})^2(\mu_1-\mu_2)\neq 0$, which is contradict with $D=0$. Hence, $\bar{\rho}_1>b_1$.
        The remainder of the proof follows the same idea in the case 1-b).
\end{enumerate}
\item In the region of $D<0$, three cases exist:
  \begin{enumerate}[a)]
    \item If $\bar{\rho}_1\geq 1$, then the conclusion follows from Proposition \ref{pro10}.
    \item If $\bar{\rho}_1< 1$ and $\bar{\rho}_1\neq \sqrt{\rho_1}$, then $\bar{\rho}_1>b_1$, the rest of the proof is similar to the case 1-b).
\item If $\bar{\rho}_1< 1$ and $\bar{\rho}_1= \sqrt{\rho_1}$, then $\bar{\rho}_1=b_1\neq \eta_1$, we have
\begin{equation*}
\lim_{\bar{\rho}_1y\rightarrow 1}\left[\frac{L^{(T)}(y)}{(1-\bar{\rho}_1y)^{-1}}\right]=L_3^{(0)}(1/\bar{\rho}_1)+\frac{\mu_1-\mu_2}{\mu_1\bar{\rho}_1}L_2^{(0)}(1/\bar{\rho}_1).
\end{equation*}
In addition, $L^{(T)}(y)$ is analytic in $\tilde{\Delta}(\phi,\varepsilon,\bar{\rho}_1)$. By applying Lemma \ref{lem3}, we get
\begin{equation*}
\pi_n^T\sim C_{t,3c}(\bar{\rho}_1)^{n}.
\end{equation*}
\end{enumerate}
\end{enumerate}\hspace*{\fill}$\Box$
\end{pf}

\section{Stochastic simulation}\label{sec6}
\begin{table}[htbp]
\centering
 \caption{\label{tab1}The ratio error of $(1-\rho)W_3$ for different loads}
 \begin{tabular}{cccccc}
  \toprule
  & $\rho=0.8$ & $\rho=0.9$ & $\rho=0.95$ & $\rho=0.975$ & $\rho=0.99$ \\
  \midrule
  \emph{E} & -95.6000  & -91.1732   & -81.4272 &  -72.7356 & -20.8699 \\
  \emph{Std} & -93.4831 & -80.9886 & -68.5507 & -39.3165 & 32.0129 \\
  \bottomrule
 \end{tabular}
\end{table}
This section tests our main results in Theorem \ref{thm1} by comparing the ratio error of the waiting times and the cdfs of the queue lengths and waiting times. The ratio error was defined in \cite{Olsen2005} by
\begin{equation*}
\hspace*{\fill}\text{Ratio error}=\frac{\text{Estimated value}-\text{Simulated value}}{\text{Simulated value}}\times 100\%, \hspace*{\fill}
\end{equation*}
where the Estimated value is the result in Theorem \ref{thm1} and the Simulated value is obtained by simulating under different traffic loads.

We consider a model with fixed parameters $\lambda_1=0.1$, $\lambda_2=0.3$, $\mu_1=0.5$, $\mu_2=1$, $\mu_3=1.5$ and $N=10$. We let $\rho=0.8, 0.9, 0.95, 0.975, 0.99$ to describe the procedure of $\rho\rightarrow 1$ and $\lambda_3$ can be determined by $\lambda_3=\mu_3(\rho-\rho_1-\rho_2)$. We use Matlab to undertake simulations under different traffic loads and each simulation runs until at least 10000 customers are served.

\begin{figure}[htp]
\centering
\includegraphics[width=4.3in]{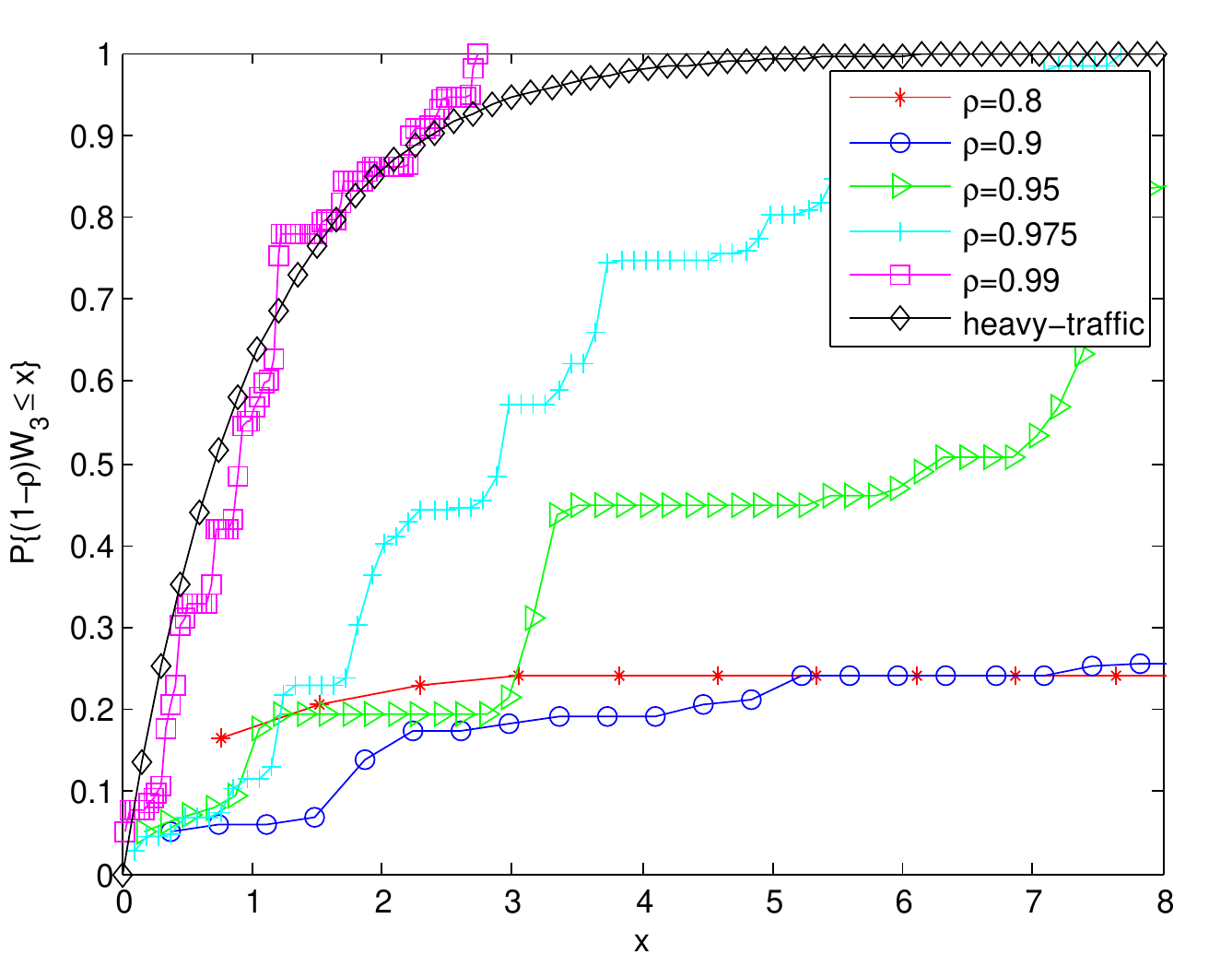} % 调用图片的时候图片的名字之间不能有任何空格也不能调用子文件夹
\caption {\small The cdf of $(1-\rho)W_3$ for different loads}
\label{fig.1}
\end{figure}

%\begin{figure}[htp]
%\centering
%\includegraphics[width=4.5in]{lengthcdf} % 调用图片的时候图片的名字之间不能有任何空格也不能调用子文件夹
%\caption {\small Empirical cumulative distribution function of $(1-\rho)X_3$ for different load}
%\label{fig.2}
%\end{figure}

For this model, the scaled queue-length and waiting-time in the critically loaded queue are exponential distributed with parameter $\eta$ and $\mu_3\eta$ respectively in the heavy-traffic scenario. Fig.\ref{fig.1} shows the cdf of $(1-\rho)W_3$ and Tab.\ref{tab1} presents the ratio error of $E(1-\rho)W_3$ and $Std(1-\rho)W_3$, where $EX$ means the expectation of $X$ and $StdX$ means the standard deviation of $X$.

It is showed that the approximation performs well when $\rho$ is very close to $1$. However when $\rho$ is moderate, the approximation seems not so accurate. This may own to the error of the simulation technique and the approximation theory since we only take the lowest order terms in the Taylor expansion. Fortunately, the higher-order terms can be obtained in the same procedure.

\begin{figure}[htp]
\centering
\includegraphics[width=5.0 in]{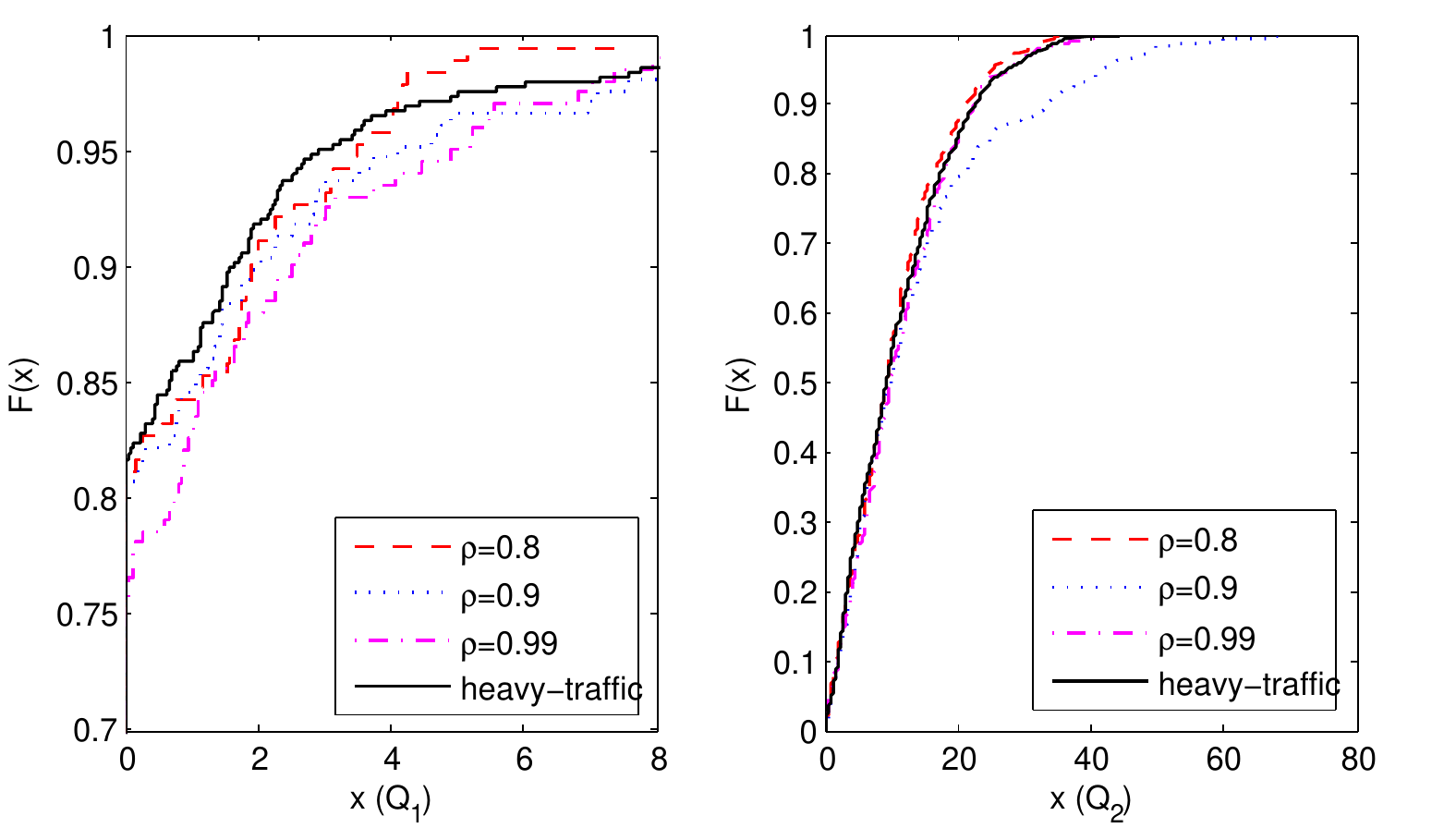} % 调用图片的时候图片的名字之间不能有任何空格也不能调用子文件夹
\caption {\small Empirical cdf of waiting times in $Q_1$ and $Q_2$ for different loads}
\label{fig.3}
\end{figure}

In the heavy-traffic regime, the queue lengths in the stable queues have the same distributions as that of a preemptive priority polling system with vacation, which is showed in Fig.\ref{fig.3}. From Fig.\ref{fig.3}, the distributions remain so closely whatever the traffic load $\rho$ is, which can be explained by the preemptive priority service policy. The queue lengths in the stable queues are independent of the value of $\rho$, which may illustrate the conclusion that the queue lengths in the stable queues and the queue length in the critically loaded queue are independent. This can be showed more exactly in non-preemptive policy systems.

 \section{Conclusions}\label{sec7}
 In this paper, we have derived the exact heavy-traffic limits of a three-queue priority polling system with threshold service policy using the singular-perturbation technique. We also provided an approximation of the tail asymptotics of the stable queues, which describes the heavy-traffic behaviors more distinctly.

 The singular-perturbation technique is based on the balance equation and hence can be extended to polling systems with multiple queues easily. It can be used to analyze the heavy-traffic limits of polling systems without multitype branching properties \cite{Resing1993}. However, if we apply the singular technique to the models with more than one critically loaded queues, then, initially, we may need to know the relative stabilities and, further, the degree of stability of each queue, which can be referred to \cite{Sum2006}. In this way, we then find the most critically loaded queue and apply the technique. In addition, when $\rho$ is moderate, the approximation seems not so accurate. Hence, it is necessary to seek for more efficient approximation techniques.

\bibliographystyle{elsarticle-num}
\bibliography{heavy-tail}

\begin{thebibliography}{10}
\expandafter\ifx\csname url\endcsname\relax
  \def\url#1{\texttt{#1}}\fi
\expandafter\ifx\csname urlprefix\endcsname\relax\def\urlprefix{URL }\fi
\expandafter\ifx\csname href\endcsname\relax
  \def\href#1#2{#2} \def\path#1{#1}\fi

\bibitem{Lee1993}
D.~S. Lee, B.~Sengupta, {Queueing analysis of a threshold based priority scheme
  for ATM networks}, IEEE/ACM Transactions on Networking (TON) 1~(6) (1993)
  709--717.

\bibitem{Boxma1995}
O.~J. Boxma, G.~M. Koole, I.~Mitrani, {A two-queue polling model with a
  threshold service policy}, in: Modeling, Analysis, and Simulation of Computer
  and Telecommunication Systems, 1995. MASCOTS'95., Proceedings of the Third
  International Workshop on, IEEE, 1995, pp. 84--88.

\bibitem{Boxma1995a}
O.~Boxma, G.~Koole, I.~Mitrani, Polling models with threshold switching, in:
  Quantitative Methods in Parallel Systems, Springer, 1995, pp. 129--140.

\bibitem{Deng2001}
Y.~Deng, J.~Tan, {Priority queueing model with changeover times and switching
  threshold}, Journal of Applied Probability 38 (2001) 263--273.

\bibitem{Deng2001a}
Y.~Deng, S.~Song, J.~Tan, {Non-Preemptive Priority queueing model with
  changeover times and switching threshold}, Communication on Applied
  Mathematics and Computation 15 (2001) 28--40.

\bibitem{Feng2001}
F.~Wei, M.~Kowada, K.~Adachi, {Performance Analysis of a Two-queue Model with
  an (M, N)-threshold Service Schedule}, Journal of the Operations Research
  Society of Japan-Keiei Kagaku 44~(2) (2001) 101--124.

\bibitem{Liu2014}
Z.~Liu, Y.~Chu, J.~Wu, {On the Three-queue Priority Polling System with
  Threshold Service Policy}, Submitted to Applied Mathematics and Computation.

\bibitem{Landry1993}
R.~Landry, I.~Stavrakakis, {Queueing study of a 3-priority policy with distinct
  service strategies}, IEEE/ACM Transactions on Networking (TON) 1~(5) (1993)
  576--589.

\bibitem{Morrison2010}
J.~A. Morrison, S.~C. Borst, {Interacting queues in heavy traffic}, Queueing
  Systems 65~(2) (2010) 135--156.

\bibitem{Boon2013}
M.~Boon, E.~Winands, Heavy-traffic analysis of k-limited polling systems, Tech.
  rep., Technical Report 2013-002, Eurandom Preprint Series, 2013. To appear in
  Probability in the Engineering and Informational Sciences. Available at
  http://www. eurandom. tue. nl/reports (2013).

\bibitem{Li2009}
H.~Li, Y.~Q. Zhao, {Exact tail asymptotics in a priority
  queue-Characterizations of the preemptive model}, Queueing Systems 63~(1-4)
  (2009) 355--381.

\bibitem{Li2011}
H.~Li, Y.~Q. Zhao, {Exact tail asymptotics in a priority queue-
  Characterizations of the non-preemptive model}, Queueing Systems 68~(2)
  (2011) 165--192.

\bibitem{Bertsimas1997}
D.~Bertsimas, G.~Mourtzinou, {Multiclass queueing systems in heavy traffic: An
  asymptotic approach based on distributional and conservation laws},
  Operations Research 45~(3) (1997) 470--487.

\bibitem{Flajolet1990}
P.~Flajolet, A.~Odlyzko, {Singularity analysis of generating functions}, SIAM
  Journal on discrete mathematics 3~(2) (1990) 216--240.

\bibitem{Olsen2005}
T.~L. Olsen, R.~D. van~der Mei, {Polling systems with periodic server routing
  in heavy traffic: renewal arrivals}, Operations Research Letters 33~(1)
  (2005) 17--25.

\bibitem{Resing1993}
J.~A.~C. Resing, {Polling systems and multitype branching processes}, Queueing
  Systems 13~(4) (1993) 409--426.

\bibitem{Sum2006}
L.~Sum, R.~K. Chang, Y.~Xie, Relative stability analysis of multiple queues,
  in: Proceedings of the 1st international conference on Performance evaluation
  methodolgies and tools, ACM, 2006, p.~65.

\end{thebibliography}
\end{document}